\documentclass[11pt]{article}
\usepackage{amsmath, amssymb, bm, url, graphics}

\newcommand{\ch}[1]{{#1}}
\newcommand{\switchText}[2]{\textcolor{red}{#1}\textcolor{blue}{#2}} 
\renewcommand{\switchText}[2]{{#2}} 

\makeatletter
\def\url@leostyle{%
  \@ifundefined{selectfont}{\def\UrlFont{\sf}}{\def\UrlFont{\small\ttfamily}}}
\makeatother
\numberwithin{equation}{section}

\newcommand{\shift}{\!\!\!\!}
\newcommand{\eps}{\varepsilon}
\newcommand{\mPoisson}[1]{\RR{\mP[#1]}}

\newcommand{\sN}{{\sss N}}
\newcommand{\qed}{\hfill $\Box$}
\newcommand{\whp}{{\bf whp}}
\newcommand{\prob}{\mathbb P}
\newcommand{\expec}{\mathbb E}
\newcommand{\eq}{\begin{equation}}
\newcommand{\en}{\end{equation}}
\newcommand{\nn}{\nonumber}

\newtheorem{theorem}{Theorem}[section]
\newtheorem{lemma}[theorem]{Lemma}
\newtheorem{prop}[theorem]{Proposition}
\newtheorem{cor}[theorem]{Corollary}
\newcommand{\proof}{\noindent {\bf Proof}. \hspace{2mm}}

\newcommand{\sss}   { \scriptscriptstyle }
\newcommand{\smallsup}[1] {{\scriptscriptstyle{(\kern-.03em{#1}}\kern-.03em)}}
\newcommand{\smallsdown}[1] {{\scriptscriptstyle{{\kern-.04em{#1}}}}}

\newcommand{\indic}[1]{{\bf 1}_{\bra{#1}}} 
\newcommand{\mindic}[1]{{\bf 1}_{#1}} 
\newcommand{\br}[1]{\left( #1 \right)}       
\newcommand{\bra}[1]{\left\{ #1 \right\}}  
\newcommand{\brc}[3]{\left#1 #3 \right#2}  
\newcommand{\brh}[1]{\left[ #1 \right]}    
\newcommand{\bfloor}[1]{\lfloor #1 \rfloor} 
\newcommand{\bceil}[1]{\lceil #1 \rceil}  
\newcommand{\mO}[1]{{\cal O}\hspace{-.0em}\bigl( #1 \bigr)} 

\newcommand{\mprob}[1]{\mathbb P\br{#1}}

\newcommand{\mexpec}[1]{\mathbb E\hspace{-.2em}\brh{#1}}

\newcommand{\mL}[1]{{\Lambda}\hspace{-.05em}_{#1}}

\newcommand{\nill}{} 

\newcommand{\ml}[1]{{\lambda}\hspace{-.05em}_{#1}}
\newcommand{\mli}[1]{\lambda \hspace{-.25em} \br{#1}}

\newcommand {\convp}{\stackrel{\sss {\mathbb P}}{\longrightarrow}}

\title{Universality for the distance in finite variance random graphs\switchText{}{:\\ Extended version}}
\author{Henri van den Esker\thanks{Delft University of Technology,
Electrical Engineering, Mathematics and Computer Science, P.O. Box 5031, 2600 GA Delft,
The Netherlands. E-mail: {\tt H.vandenEsker@ewi.tudelft.nl} and {\tt G.Hooghiemstra@ewi.tudelft.nl}}
\and
Remco van der Hofstad
\thanks{Department of Mathematics and
Computer Science, Eindhoven University of Technology, P.O.\ Box
513, 5600 MB Eindhoven, The Netherlands. E-mail: {\tt
rhofstad@win.tue.nl}}
\and Gerard Hooghiemstra$^*$}

\begin{document}
    \maketitle

\begin{abstract}
We generalize the asymptotic behavior of the graph distance
between two uniformly chosen nodes in the configuration model to a
wide class of random graphs. Among others, this class contains the
Poissonian random graph, the expected degree random graph and the
generalized random graph (including the classical Erd\H{o}s-R\'enyi graph).

In the paper we assign to each node a deterministic capacity and
the probability that there exists an edge between a pair of nodes
is equal to a function of  the product of the capacities of the
pair divided by the total capacity of all the nodes. \ch{We
consider capacities which are such that the degrees of a node has
uniformly bounded moments of order strictly larger than two, so
that, in particular, the degrees have finite variance.} We prove
that the graph distance grows like $\log_\nu N$, where the $\nu$
depends on the capacities. In addition, the random fluctuations
around this asymptotic mean $\log_\nu N$ are shown to be tight. We
also consider the case where the capacities are independent copies
of  a positive random $\mL\nill$ with $\mprob{\mL\nill > x}\leq
cx^{1-\tau}$, for some constant $c$ and $\tau>3$, again resulting
in graphs where the degrees have finite variance.

The method of proof of these results is to couple each member
of the class to the Poissonian random graph, for which we then give
the complete proof by adapting the arguments of \cite{finstub}.
\end{abstract}

\newcommand{\Q}{\mathbb Q}
\newcommand{\mCo}[1]{{\bf C#1}}
\newcommand{\mRNprocess}{$\overline{\mathrm{NR}}$-process} 
\newcommand{\mrRNprocess}{$\underline{\mathrm{NR}}$-process} 

\newcommand{\mP}[1][\mL\nill]{Poi\hspace{-.15em}\br{#1}}

\newcommand{\SPT}{{\cal N}}

\newcommand{\eqn}[1]{\begin{equation} #1 \end{equation}}
\newcommand{\eqan}[1]{\begin{align} #1 \end{align}}
\newcommand{\vep}{\varepsilon}

\newcommand{\mh}[1]{t_{#1}\hspace{-.1em}}

\newcommand{\msM}[2][k]{{\cal M}_{\smallsdown{#1,#2}}} 
\newcommand{\mM}[2][k]{\overline M_{\smallsdown{#1,#2}}} 
\newcommand{\mMi}[3][k]{\overline M^{\smallsup{#2}}_{\smallsdown{#1,#3}}} 
\newcommand{\mnM}[2][k]{{\underline M}_{\smallsdown{#1,#2}}} 
\newcommand{\mnMi}[3][k]{{\underline M}^{\smallsup{#2}}_{\smallsdown{#1,#3}}} 
\newcommand{\mMd}[3]{M^{\smallsup{#1}}_{\smallsdown{#2,#3}}}


\newcommand{\RR}[1]{\textcolor{red}{#1}} 
\renewcommand{\RR}[1]{#1}
\newcommand{\mf}[2][]{\RR{#1f^{\smallsup{N}}_{#2}}} 
\newcommand{\mg}[2][]{\RR{#1g^{\smallsup{N}}_{#2}}} 
\newcommand{\mfl}[1]{\RR{f_{#1}}} 
\newcommand{\mgl}[1]{\RR{g_{#1}}} 

\newcommand{\mnZ}[1]{{\underline Z}^{}_{\smallsdown{#1}}}
\newcommand{\mnZi}[2]{{\underline Z}_{\smallsdown{#2}}^{\smallsup{#1}}}
\newcommand{\mZ}[1]{{\overline Z}_{\smallsdown{#1}}}
\newcommand{\mZi}[2]{{\overline Z}_{\smallsdown{#2}}^{\smallsup{#1}}}
\newcommand{\mZd}[2]{{Z}_{\smallsdown{#2}}^{\smallsup{#1}}}

\newcommand{\msZ}[1]{{\cal Z}_{\smallsdown{#1}}}
\newcommand{\msZi}[2]{{\cal Z}_{\smallsdown{#2}}^{\smallsup{#1}}}

\newcommand{\mC}[1]{{\overline C}_{\smallsdown{#1}}}
\newcommand{\mCi}[2]{{\overline C}_{\smallsdown{#2}}^{\smallsup{#1}}}

\newcommand{\mnC}[1]{{\underline C}_{\smallsdown{#1}}}
\newcommand{\mnCi}[2]{{\underline C}_{\smallsdown{#2}}^{\smallsup{#1}}}

\newcommand{\mJ}[1]{ {{\overline {\bm M}}}_{\smallsdown{#1}}}
\newcommand{\mJd}[2]{{ { {\bm M}}}^{\smallsup{#1}}_{\smallsdown{#2}}}
\newcommand{\mJi}[2]{{\overline {\bm M}}_{\smallsdown{#2}}^{\smallsup{#1}}}
\newcommand{\mnJ}[1]{{ { \underline {\bm M}}}_{\smallsdown{#1}}}
\newcommand{\mnJi}[2]{{\underline {\bm M}}_{\smallsdown{#2}}^{\smallsup{#1}}}

\section{Introduction}
Various papers (see e.g.,
\cite{BBCR,CLaverage,finstub,NSW00,norros3}) study properties of
random graphs with a given degree sequence. Among such properties
as connectivity, cluster size and diameter, the graph distance
between two uniformly chosen nodes is an important one. For two
connected nodes the graph distance is defined as the minimum
number of edges of a path that connects these nodes. If the nodes
are not connected, then the graph distance is put equal to
infinity.

For the configuration model (see Section 1.4 for a definition) a
distance result appeared in \cite{finstub}, when the distribution
of the i.i.d.~degrees $D^{\smallsup{C}}$ satisfies
\begin{equation}
\label{degreedistr}
\prob(D^{\smallsup{C}}>x)\leq c x^{1-\tau},
\end{equation}
for some constant $c$, all $x\geq 0$, and with $\tau>3$.
We use the superscript $(C)$ to differentiate between models.
The result in \cite{finstub} states that with probability converging to $1$ ({\bf whp}), the \ch{average} distance between nodes in the giant
component  has, for
\begin{equation}
\label{def-nu}
\tilde \nu=\frac{\expec[D^{\smallsup{C}}(D^{\smallsup{C}}-1)]}{\expec[D^{\smallsup{C}}]}>1,
\end{equation}
 bounded fluctuations around $\log_{\tilde \nu}N$. The condition $\tilde \nu>1$ corresponds to the supercritical
 case of an associated branching process, and is the condition
 under which a giant component exists.

In this paper we extend the above distance result to a wide class of random graphs.
Models which fall in this class  are the generalized random graph (GRG), the expected degree random graph
(EDRG) and the Poissonian random graph (PRG). All three models will be introduced in more detail below.

The method of proof is coupling. It is shown that the distance
result holds for all models in the general class if and only if
the result holds for the PRG (Section 2). In Section 4 we prove
the distance result for the Poissonian random graph. This proof is
parallel to that in \cite{finstub} for the configuration model.
\switchText{Therefore, most of the proofs of the auxiliary propositions and lemmas are left out, as they are similar to those in \cite{finstub}. Since details of these proofs are different, we included full proofs in an extended version of this paper. The extended version will not be published, but is available on the web  (\cite{EV}).}
    {In this paper we included full proofs of the auxiliary lemmas contained in Section \ref{NRmodel},
    since details of these proofs are different from those in \cite{finstub}.}

\subsection{Model assumptions} \label{sec:model definition}
The graph models considered here are \emph{static} models, meaning
that the number of nodes is fixed. The graph $G_\sN$ has $N$
nodes, numbered $1,2,\ldots,N$. Associated with the nodes is a
sequence $\{\ml i\}_{i=1}^N$ of positive reals. We call $\ml i$
the \emph{capacity} of node~$i$; nodes with  a large capacity will obtain a high degree, whereas nodes with small capacity
will only be incident to a limited number of edges.
Furthermore, we define
\begin{align} \label{LsN}
l_\sN =\ml 1+\ml 2+\cdots+\ml \sN,\,
\end{align}
i.e., $l_\sN$ is the total capacity of all nodes of the graph $G_\sN$.

The binary random variables $\{X_{ij}\}_{1\leq i \leq j \leq N}$, are defined by setting $X_{ij}=1$,
if the edge between node $i$ and node $j$ in the graph $G_{\sN}$ is present, otherwise we set $X_{ij}=0$.
If $ i> j$, then  by convention $X_{ji}=X_{ij}$.
We call $\{X_{ij}\}$ the \emph{connection variables} and $p_{ij}=\mprob{X_{ij}=1}$ the \emph{connection probability
of edge $ij$}.
In this paper we consider graphs $G_\sN$ that satisfy the following two assumptions:
     \begin{enumerate}
            \item[{\bf A1}:] The connection variables $\{X_{ij}\}_{1\leq i < j \leq N}$, are independent.

        \item[{\bf A2}:] The connection probability $p_{ij}$, for $1\leq i < j \leq N$, can be written as
        $p_{ij}=h({\ml i\ml j}/{l_\sN})$, for some function $h:[0,\infty)\rightarrow[0,1]$, satisfying
        \begin{align}\label{conditions''}
        h(x)-x =\mO{x^{2}},& \,\,\,\,\text{for } x \downarrow 0.
        \end{align}
     \end{enumerate}

\subsection{Three special cases} \label{subseq:sepc cases}
We give three examples of random graph models, which satisfy assumptions {\bf A1} and {\bf A2}, and
hence fall in the class of models considered here.

The first example is the Poissonian random graph (PRG), which is a variant of a model introduced by Norros and Reittu in \cite{norros3}.
The second  random graph model, which we call the expected degree random graph (EDRG),
is introduced by Chung and Lu in \cite{CLaverage,CLconnected}.
The third and last example is the generalized random graph (GRG), which was introduced by Britton, Deijfen and Martin-L\"of in \cite{Britton}.

We now define the three models and verify that each of them satisfy the conditions {\bf A1} and {\bf A2}.
\begin{itemize}
        \item{\bf Poissonian random graph:}
        In \cite{norros3} the Poissonian random graph (PRG) was introduced.
        The main feature of such a graph $G^{\smallsup{P}}_\sN$ is that,
        conditionally on the capacities, the number of edges between any pair of nodes $i$ and $j$ is a Poisson random variable.
        The model in \cite{norros3} is introduced as a growth model, but as a consequence of  \cite[Proposition 2.1]{norros3},
        it can be formulated as a static model, and we will do so. Start with the graph $G^{\smallsup{P}}_\sN$ consisting of $N$ nodes
        and capacities $\{\ml i\}_{i=1}^N$.
        The number of edges between two different nodes $i$ and $j$ is given by an independent Poisson random variable
        $E^{\smallsup{P}}_{ij}$ with parameter
        \begin{align}\label{poissonedge}
        \ml i\ml j/l_\sN.
        \end{align}

      Denote by $\mindic{A}$ the indicator of the set $A$. The connection variables are then $X^{\smallsup{P}}_{ij}=\mindic{\{E^{\smallsup{P}}_{ij}>0\}}$, so that, for $1 \leq i \leq j  \leq N$,
      the connection probabilities are given by
     $$
     p^{\smallsup{P}}_{ij}=\prob(X^{\smallsup{P}}_{ij}=1)=\prob(E^{\smallsup{P}}_{ij}>0)=
     1-\exp(-\ml i\ml j/l_\sN)=h^{\smallsup{P}}({\ml i\ml j}/{l_\sN}),
     $$
      where $h^{\smallsup{P}}(x)=1-e^{-x}$. Obviously, $h^{\smallsup{P}}(x)-x =\mO{x^2},$ for $x\downarrow 0$. Since, by definition,
      the random variables $\{X^{\smallsup{P}}_{ij}\}_{1\leq i < j \leq N}$ are independent,
       we conclude that the assumptions {\bf A1} and {\bf A2} are satisfied.

        It should be noted that in this paper we define the PRG using a deterministic sequence of capacities, while Norros and Reittu start
        with an i.i.d.~sequence of random capacities. The case where the capacities $\{\Lambda_i\}_{i=1}^N$ are i.i.d.~random
        variables, satisfying certain tail estimates, is a special case of our set up and is studied in more detail in Theorem \ref{corIIDcase}.

                \item {\bf Expected degree random graph:}
                In \cite{CLaverage,CLconnected} a random graph model is introduced
                starting from a sequence of deterministic capacities $\{\ml i\}_{i=1}^N$.
                We construct the EDRG $G_\sN^{\smallsup{E}}$ as follows. Let $\{X^{\smallsup{E}}_{ij}\}_{1 \leq i \leq j \leq N}$  be a sequence of independent Bernoulli random variables with
                success probability
                $$
                p^{\smallsup{E}}_{ij}=\mprob{X^{\smallsup{E}}_{ij}=1}=(\ml i\ml j / l_\sN)\wedge 1, \,\, \text{ for }1 \leq i \leq j \leq N,
                $$
                where $x\wedge y$ denotes the minimum of $x$ and $y$.
                This minimum is to ensure that the result is a probability.

        Assumption {\bf A1} is satisfied by definition, since the connection variables  are independent Bernoulli variables,
        and assumption {\bf A2} is also satisfied if we pick $h^{\smallsup{E}}(x)=x \wedge 1$.

                If we assume that $\ml i\ml j/l_\sN<1$ for all $1 \leq i \leq j \leq N$,
                then the expected degree of a node $i$ is given by $\ml i$:
        $$
        \mexpec{\sum_{j=1}^N X^{\smallsup{E}}_{ij}}=\sum_{j=1}^N \ml i\ml j/l_\sN=\ml i.
        $$

    The Erd\H{o}s-R\'enyi random graph, usually denoted by $G(N,p)$, is a special case of the EDRG.
    In the graph $G(N,p)$, an edge between a pair of nodes is present with probability $p\in[0,1]$, independently of the other edges.
    When $p=\lambda/N$, for some constant $\lambda>0$,
    then we obtain the graph $G(N,\lambda/N)$ from the EDRG by picking $\ml i=\lambda$, for all $i$, since then
    $p_{ij}^{\smallsup{E}}=\ml i \ml j/ l_\sN=\lambda/N=p$, for all $1 \leq i < j \leq N$.

    \item {\bf Generalized random graph:}
    The generalized random graph (GRG) is an adapted version of the EDRG, see the previous example.
    We define $G_\sN^{\smallsup{G}}$ with $N$ nodes as follows.
    The sequence of connection variables is again given by a sequence of independent Bernoulli
    random variables $\{X^{\smallsup{G}}_{ij}\}_{1\leq i < j \leq N}$ with
    $$
    \mprob{X^{\smallsup{G}}_{ij}=1}=p^{\smallsup{G}}_{ij}=\frac{\ml i \ml j/l_\sN}{1+\ml i \ml j/l_\sN}.
    $$
    In \cite{Britton} the edge probabilities are given by $(\ml i\ml j/N)/(1+\ml i\ml j/N)$,
    so that we have replaced $\ml i/N^{1/2}$   by $\ml i/l_\sN ^{1/2}$, $1\le i \le N$. This makes hardly any difference.

Again, the assumptions {\bf A1} and {\bf A2} are satisfied.
To satisfy assumption {\bf A2} we pick $h^{\smallsup{G}}(x)=x/(1+x)=x+\mO{x^2}$.
\end{itemize}

\subsection{Main results} \label{sec: main results}
We state conditions on the capacities $\{ \ml i\}_{i=1}^N$, under which our main result will hold. We shall need three conditions,
which we denote by (\mCo1), (\mCo2) and (\mCo3), respectively.

\paragraph{(C1) Convergence of means:} Define
    \begin{align}
    \label{muNnuN-def}
    \mu_{\sN}&=\frac 1N \sum_{i=1}^N \lambda_i&
    &\text{and}&
    \nu_{\sN}&=\frac{\sum_{i=1}^N \lambda^2_i}{\sum_{i=1}^N \lambda_i},&
    \end{align}
 then there exist constants $\mu\in (0,\infty)$, $\nu\in (1,\infty)$ and $\alpha_1>0$ such that
    \eqn{
    \label{reqA}
    |\mu_{\sN}-\mu|=\mO{N^{-\alpha_1}}, \qquad |\nu_{\sN}-\nu|=\mO{ N^{-\alpha_1}}.
    }
\paragraph{(C2) Convergence of branching processes:} Define
    \begin{align}
    \label{def-fgn-n}
    \mf n
    &=\frac{1}{N}\sum_{i=1}^N e^{-\lambda_i}\frac{\lambda_i^n}{n!}&
    &\text{and}&
    \mg n&=\frac{1}{N\mu_{\sN}}\sum_{i=1}^N e^{-\lambda_i}\frac{\lambda_i^{n+1}}{n!},&
    \end{align}
then there exist sequences $\{\mfl n\}_{n\geq 0}$ and $\{\mgl n\}_{n\geq 0}$, independent of $N$,
and $\alpha_2>0$ such that
    \begin{align}
    \label{reqB}
   {\rm d}_{\sss \rm TV}(\mf\nill, \mfl\nill) &=\mO{N^{-\alpha_2}}&
    &\text{and}&
   {\rm d}_{\sss \rm TV}(\mg\nill, \mgl\nill) &=\mO{N^{-\alpha_2}},&
    \end{align}
    where ${\rm d}_{\sss \rm TV}(\,\cdot\,,\,\cdot\,)$ is the total variance distance, i.e., for probability mass functions $p=\{p_j\}$
    and $q=\{q_j\}$:
    \begin{equation}
    \label{vardistance}
{\rm d}_{\sss \rm TV}(p,q)=\tfrac12 \sum_j |p_j-q_j|.
\end{equation}

\paragraph{(C3) Moment and maximal bound on capacities:}
There exists a $\tau>3$ such that
for every $\vep>0$,
    \begin{align}
    \label{reqC}
    &\limsup_{N\rightarrow \infty} \frac{1}{N}\sum_{i=1}^N \ml i^{\tau-1-\varepsilon}<\infty,&
    &\text{and}&
    \lambda_\sN^\smallsup N&=\max_{1 \leq i \leq N} \lambda_i \leq N^{\gamma},&
    \end{align}
    where,
    \eqn{\label{P:assumP}
    \gamma=\frac1{\tau-1}+\vep < 1/2.
    }
It is not hard to see that $\mu_{\sN}$ and $\nu_{\sN}$ in
\eqref{muNnuN-def} are the means of the probability mass functions
$\{\mf n\}_{n\geq0}$, $\{\mg n\}_{n\geq0}$,
respectively. Thus, (\mCo1) is equivalent to the fact that the means
of $\{\mf n\}_{n\geq0}$, $\{\mg n\}_{n\geq0}$
converge. \ch{It turns out that (\mCo1) is equivalent to
the convergence of the first and second moment of the degree
of a uniform node.}

{\ch{Condition} (\mCo2) says that the laws $\{\mf n\}_{n\geq0}$,
$\{\mg n\}_{n\geq0}$ are close to certain limiting laws
$\{\mfl n\}_{n\geq0}$, $\{\mgl n\}_{n\geq0}$, \ch{which shall turn out to be
crucial in our proofs, since it allows us to use coupling to branching processes.}

\ch{The second bound in Condition (\mCo3) gives an upper bound
on the maximal capacity of a node, while it can be seen that
the first inequality is equivalent to the statement that a uniform
node has a uniformly bounded moment of order at least $\tau-1-\varepsilon$.
Since $\tau>3$, we can pick $\varepsilon>0$ so small that $\tau-1-\varepsilon>2$,
so that the degrees have {\it finite variances}.}

\ch{We shall prove our main results in the generality of Conditions (\mCo1)--(\mCo3),
but shall give special cases when Conditions (\mCo1)--(\mCo3) following our main results.}

\ch{In order to be able to state our main results, we}
define the process $\{\msZ t\}_{t\geq0}$ as a \ch{branching process (BP)}
starting from $\msZ 0=1$, where in the first generation
the offspring distribution is equal to $\{\mfl n\}$, whereas
in the second and further generations the offspring is chosen in
accordance to $\{\mgl n\}$.

We define the \ch{average} graph distance or hopcount
$H_\sN$ between two different randomly chosen nodes $A_1$
and $A_2$ in the graph $G_\sN$ as the minimum number of edges
that form a path from the node $A_1$ to node $A_2$ where, by
convention, the distance equals $\infty$ if the nodes $A_1$
and $A_2$ are not connected.

\begin{theorem}[Fluctuations of the graph distance]
\label{main}\label{P:UNIV:main}
Assume that the capacities $\{\ml i\}_{i=1}^N$
satisfy (\mCo1)--(\mCo3), and let the graph $G_N$, with
capacities  $\{\ml i\}_{i=1}^N$ satisfy  {\bf A1} and {\bf A2}, for some function
$h: [0,\infty)\mapsto [0,1]$.
Let $\sigma_N=\bfloor{\log_\nu N}$ and $a_N=\sigma_N -\log_\nu N$. There exist random variables $(R_a)_{a\in(-1,0]}$ such that, as
$N\rightarrow\infty$,
\begin{align}
\label{annealed}
\mprob{H_\sN=\sigma_\sN+k \, | \, H_\sN < \infty}=\mprob{R_{a_{\sN}}=k}+o(1), \qquad k=0,\pm1,\pm2,\ldots.
\end{align}
\end{theorem}

We identify the random variables $(R_a)_{a\in(-1,0]}$ in Theorem \ref{limit law} below. Before doing so,
we state one of the consequences of Theorem \ref{main}:
\begin{cor}[Concentration of the graph distance] \label{tightness} Under the given assumptions of Theorem \ref{main},
\begin{itemize}
    \item with probability $1-o(1)$ and conditionally on $H_\sN<\infty$, the random variable $H_\sN$ is in between $(1\pm\eps)\log_\nu N$ for any $\eps>0$;

    \item conditionally on $H_\sN<\infty$, the sequence of random variables $H_\sN-\log_\nu N$ is tight, i.e.,
    \eq
    \lim_{K\rightarrow \infty}\limsup_{N\rightarrow \infty}\mprob{|H_\sN-\log_\nu N| \leq K \bigl.\bigr| H_\sN <\infty}=1.
    \en
\end{itemize}
\end{cor}

We use a limit result from branching process theory to identify the limiting random variables $(R_a)_{a\in(-1,0]}$.
It is well known, see \cite[p.~244]{fellerb}, that the process $\{\msZ t/\mu\nu^{t-1}\}_{t\geq1}$ is a
non-negative martingale and consequently converges almost surely to a limit ${\cal W}$:
\begin{align}  \label{lawW}
\lim_{t\rightarrow\infty} \frac{\msZ t}{\mu\nu^{t-1}}={\cal W}, \text{ \rm a.s.}
\end{align}

Let ${\cal W}^{\smallsup{1}}$ and ${\cal W}^{\smallsup{2}}$ be two independent copies of
${\cal W}$ in \eqref{lawW}, then we can identify the limit random variables $(R_{\smallsdown{a}})_{a\in(-1,0]}$ as follows:
\begin{theorem}
\label{limit law}
    For $a\in(-1,0]$,
    $$
    \mprob{R_{\smallsdown{a}} > j}=\mexpec{\exp\bra{-\kappa \nu^{a+j}{\cal W}^{\smallsup{1}}{\cal W}^{\smallsup{2}}}|{\cal W}^{\smallsup{1}}{\cal W}^{\smallsup{2}}>0},
    $$
    where $\kappa=\mu(\nu-1)^{-1}$.
\end{theorem}

A special case of the above theorems is the case where
\begin{equation}
\label{speccase}
\lambda_i=(1-F)^{-1}(i/N),\qquad i=1,2,\ldots,N,
\end{equation}
where $F$ is a distribution function of a positive random variable, i.e.,
$F(0)=0$, $F$ is non-decreasing and $\lim _{x\to\infty } F(x)=1$, see \eqref{invverd} for a definition of $(1-F)^{-1}$.
In Appendix B, we will formulate quite general conditions on a distribution function $F$
such that (\mCo1)--(\mCo3) hold with $\{\ml i\}_{i=1}^N$ defined by \eqref{speccase}. The special cases
$$
F(x)=1-\frac{c}{x^{\tau-1}},\qquad \mbox{for which} \qquad (1-F)^{-1}(y)=\Big(\frac{c}{y}\Big)^{1/(\tau-1)},
$$
where $\tau >3$, \ch{extends the results obtained by Chung and Lu \cite{CLaverage}
from $H_{\sN}/\log_\nu{N}\convp 1$ to the study of the
fluctuations of $H_{\sN}$.}

Instead of assigning {\it deterministic} capacities to the nodes
one can also assign {\it random} capacities.
Associated with the nodes is a sequence $\{\mL i\}_{i=1}^N$
of positive i.i.d. random variables, with distribution
    $$
    F_\Lambda(x)=\mprob{\mL\nill \leq x}.
    $$
Then, we set, for $1 \leq i \leq N$,
    $$
    \lambda_i=\mL i.
    $$
For the i.i.d.~case, we \ch{can identify $\mu, \nu, f$ and $g$ appearing in
conditions (\mCo1)--(\mCo3) as}
    \begin{align}\label{pick idd}
        \mu& =\mexpec{\mL\nill},&
        \nu& =\frac{\mexpec{\mL\nill^2}}{\mexpec{\mL\nill}},&
    \mfl n&=\mexpec{e^{-\mL \nill}\frac{\mL \nill^n}{n!}}&
    &\text{and}&
    \mgl n&=\frac{1}{\mu}\mexpec{e^{-\mL \nill}\frac{\mL \nill^{n+1}}{n!}},&
    \end{align}
    for $n\geq 0$.

The next theorem states that results of the deterministic
capacities carry over to the random capacities:
\begin{theorem} \label{corIIDcase}
    Given an i.i.d.~sequence of random variables $\{\mL i\}_{i=1}^N$,
   with common distribution function $F_\Lambda$.
 If there exist constants $c>0$ and $\tau>3$ such that
\begin{align}
\label{distribution} &1-F_\Lambda(x) \leq c x^{1-\tau} \text{,
for all } x \geq 0,
\end{align}
and with $\nu$, given by \eqref{pick idd}, satisfying $\nu>1$,
then there exists an event ${\cal J}$, which occurs with high probability, such that conditionally on
$\{\mL i\}_{i=1}^N$ satisfying ${\cal J}$, the conditions (\mCo1)--(\mCo3) hold.
\end{theorem}
The proof of Theorem \ref{corIIDcase} is given in Appendix
\ref{iidcase}. We have \ch{the following corollary to Theorem  \ref{corIIDcase}:}

\begin{cor}
\label{IIDcase} In the case of i.i.d.~capacities, with common
distribution function $F_{\Lambda}$ satisfying \eqref{distribution}
and with $\nu>1$, the results of Theorem \ref{P:UNIV:main}, Corollary \ref{tightness} and
Theorem \ref{limit law} hold \ch{with high probability. More precisely,
for every $k=0,\pm1,\pm2,\ldots$, the random variable
    \eqn{
    \label{quenched}
    \frac{\mprob{H_\sN=\sigma_\sN+k\,\big|\,\{\mL i\}_{i=1}^N}}{\mprob{H_\sN < \infty\,\big|\,\{\mL i\}_{i=1}^N}}-\mprob{R_{a_{\sN}}=k}
    }
converges in probability to zero.
}
\end{cor}
\vskip0.5cm

\begin{figure}
	\begin{center}
			\includegraphics{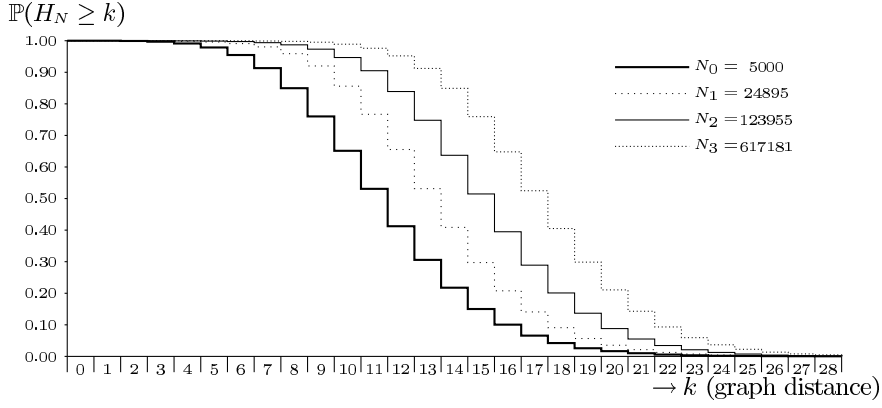}
		\end{center}
			\caption{Empirical survival functions of the graph distance for $\tau=3.5$ and for four values of $N$. Each plot is averaged over 1000 samples.\label{fig}}
\end{figure}

We demonstrate Theorem \ref{P:UNIV:main} for the i.i.d. case using Corollary \ref{IIDcase}. Assume that
$$
1-F_{\Lambda}(x)=cx^{1-\tau}\indic{x > x_0},
$$
with $\tau=3.5$, $c=2.5981847$ and $x_0=0.7437937$, then $\nu=\mexpec{\Lambda^2}/\mexpec{\Lambda}\approx 2.231381$.
We can pick different values of the size of the simulated graph, so that for each two simulated values $N$ and $M$ we have $a_N = a_M$. As an example, we take $N=M\nu^{2k}$, for some integer $k$. This induces, starting from $N_0=M=5000$, by taking for $k$ the successive values 0,1,2,3
\begin{align}
 N_0=M&=5000,&
 N_1&=24895,&
 N_2&=123955,&
 &\text{and}&
 N_3&=617181.&
\end{align}
Observe that $a_{N_k} = 0.6117...$ for $k=0,1,2,3$.
According to Corollary \ref{IIDcase}, the survival function of the graph distance $H_{\scriptscriptstyle N}$, run approximately parallel on distance 2 in the limit for $N\rightarrow \infty$, since $\log_\nu N_k=\log_\nu M + 2k$ for $i=0,1,2,3$.
In Figure \ref{fig} we have simulated the survival function  of the graph distance for the graphs with sizes $N_k$ with $k=0,1,2,3$, and, indeed, the plots are approximately parallel on distance 2.

\vskip0.5cm

\noindent
A general version of inhomogeneous random graphs with independent
edges is presented in \cite{phasetrans}. Our choice of graphs
is a special case of the rank-1 case treated in
\cite[Section 16.4]{phasetrans}. In the general setting of
\cite{phasetrans}, the vertices are given by $\{x_i\}_{i=1}^N\subset {\cal X}$,
for some state space ${\cal X}$, and there is an edge between vertices
$x_i$ and $x_j$ with probability
    \eqn{
    \label{pij-BolJanRio}
    p_{ij} = \min\{\kappa(x_i,x_j)/N,1\},
    }
where $\kappa\colon {\cal X} \times {\cal X}\to [0,\infty)$ is a suitable kernel.
The rank-1 case of \cite{phasetrans} is the case where $\kappa$ is of product form,
i.e., $\kappa(x_i,x_j)=\psi(x_i)\psi(x_j)$ for some
$\psi\colon {\cal X}\to [0,\infty)$. In fact, it is even possible
that $\kappa$ depends on $N$, i.e., $\kappa=\kappa_{\sN}$,
such that $\kappa_{\sN}$ converges to some limit as $N\rightarrow \infty$
in a suitable way. This allows one to simultaneously treat the cases where
$p_{ij}=1-e^{-\kappa(x_i,x_j)/N}$ or $p_{ij} = \kappa(x_i,x_j)/(N+\kappa(x_i,x_j))$
(recall the special cases in Section \ref{subseq:sepc cases}).
In \cite{phasetrans}, various results are proved in this generality,
such as the phase transition of the resulting graphs, and distance results,
such as \ch{average} distances and the diameter. The main tool
is a deep comparison to multitype Poisson branching processes.

In particular, \cite[Theorem 3.14]{phasetrans} states that for suitable
{\it bounded} $\kappa_{\sN}$, the \ch{average} distances between uniform pairs of
connected nodes is equal to $(1+o(1))\log_{\nu}{N}$, for a certain
$\nu>1$. The condition $\nu>1$ corresponds exactly to the
random graph being supercritical (see \cite[Theorem 3.1]{phasetrans}).
This \ch{average} distance result generalizes the first result
in Corollary \ref{tightness}, apart from the fact that
in our setting we do not assume that $\kappa_{\sN}$ is bounded.
In fact, in certain cases, $\kappa_{\sN}$ can be unbounded.
In our paper, we state conditions on $\psi$ in the rank-1
case of \cite{phasetrans} under which we
can identify the {\it fluctuations} of the \ch{average} distances.
It would be interesting to generalize our main results
to the general setting of \cite{phasetrans}, i.e.,
to study the fluctuations of $H_{\sN}$ in the general setting of
\cite{phasetrans}. However, it is unclear to us how we should
generalize the tail assumption in \eqref{distribution} to this
inhomogeneous setting.

\subsection{Relations with the configuration model}\label{P:UNIV:CompareCMandIRG}
The configuration model (CM) appeared in the context of random regular graphs as
early as 1978 (see \cite{BeC,JLR}).
Molloy and Reed \cite{MR95} were the first to use the configuration model with specified degree sequences.
Here we consider the CM as introduced in \cite{finstub}.
Start with an i.i.d. sequence $\{D_i^{\smallsup{C}}\}_{i=1}^N$ of positive integer valued
random variables , where $D_i^{\smallsup{C}}$ will denote the degree of node $i$.
To built a graph it is mandatory that $D_1^{\smallsup{C}}+D_2^{\smallsup{C}}+\ldots+D_\sN^{\smallsup{C}}$
is even, so if $D_1^{\smallsup{C}}+D_2^{\smallsup{C}}+\ldots+D_\sN^{\smallsup{C}}$ is odd we increase
$D^{\smallsup{C}}_\sN$ by one, which will have little effect.
We build the graph model by attaching $D^{\smallsup{C}}_i$ stubs or half edges to node $i$
and pair the stubs at random, so that two half edges will form one edge.

In \cite{finstub}, the authors prove a version of Theorem \ref{main}-\ref{limit law}
for the configuration model. Theorems \ref{main}-\ref{limit law} hold for
the configuration model with only two changes:
\begin{enumerate}
\item[1.] Replace the condition $\nu>1$ in Theorem \ref{main}, which is implicitly contained in (\mCo1) by the condition $\tilde \nu>1$, defined in (\ref{def-nu}).

\item[2.] Replace the offspring distributions of the BP $\{{\cal Z}_t\}_{t\ge 0}$, by
\begin{enumerate}
\item
    $$
    \tilde f_n=
    \prob(D^{\smallsup{C}}=n), \quad n\ge 1,
    $$
\item
    $$
    \tilde  g_n=\frac{(n+1)\tilde  f_{n+1}}
    {\expec{[D^{\smallsup{C}}}]}, \quad n\geq 0.
    $$
\end{enumerate}
\end{enumerate}
\ch{For the configuration model, the setting is as in Theorem \ref{corIIDcase},
where, for the CM, the {\it degrees} are chosen in an i.d.d.\ fashion.
The result in \cite{finstub} proves that when we pick two uniform nodes,
that their distance, when averaged out over the randomness in
the i.i.d.\ degrees, satisfies \eqref{annealed}.
A version in \eqref{quenched}, which holds conditionally on
the random degrees, is {\it not} proved in \cite{finstub}, and
is stronger than \eqref{annealed}.
We conjecture that a version of \eqref{quenched} also holds in
the CM, when the degrees satisfy appropriate conditions that would take
the place of (\mCo1)--(\mCo3).}

One wonders why a result like Theorems 1.1-1.3 holds true for
the class of models introduced in Section 1.1,
especially if one realizes that in the CM the degrees are
\emph{independent}, and the edges are not, whereas
for instance in the GRG (and in the other two examples)
precisely the opposite is true, i.e., in the GRG the edges are
independent and the degrees are not. To understand
at least at an intuitive level why the distance result
holds true, we compare the configuration model with the
generalized random graph.

By construction the degree sequence $D_1^{\smallsup{C}},D_2^{\smallsup{C}},\ldots,D_\sN^{\smallsup{C}}$ of the CM is an
i.i.d. sequence, and conditionally on ${\cal D}=\{D_1^{\smallsup{C}}=d_1,D_2^{\smallsup{C}}=d_2,\ldots,D_\sN^{\smallsup{C}}=d_\sN\}$,
the graph configuration is uniform over all configurations satisfying ${\cal D}$,
because the pairing is at random. Hence if we condition on both the event ${\cal D}$
and the event ${\cal S}=\{\mbox{the resulting graph has no self-loops and no multiple edges}\}$, then the CM
renders a simple graph, which is picked uniformly from all possible {\it simple} configurations
with degree sequence satisfying ${\cal D}$. Since for $N\to \infty$ the probability
of the event ${\cal S}$ converges to $\exp(-{\tilde \nu}/2-{\tilde \nu}^2/4)>0$ (see \cite[p. 51]{bollobasboek}), it follows
from \cite[Theorem 9.9]{JLR} that properties that hold {\bf whp} in the CM also
hold {\bf whp} in the conditioned simple graph. Hence a property as tightness of the
graph distance $H_\sN^{\smallsup{C}}$ in the CM is inherited by the conditioned simple graph,
with the same degree sequence.
This suggests that also the limiting distribution of the fluctuations of the graph distance in the CM conditioned on ${\cal S}$
is the same as the one in the CM as identified in \cite{finstub}. A direct proof of this claim is missing though.

On the other hand the GRG with given degree sequence $d_1,d_2,\ldots,d_\sN$ is also
uniform over all possible (simple) configurations. Moreover \cite[ Theorem 3.1]{Britton} shows that
the degree sequence ${D_1^{\smallsup{G}}},{D_2^{\smallsup{G}}},\ldots,{D^{\smallsup{G}}_\sN}$ of the GRG  is asymptotically
independent with marginal distribution a mixed Poisson distribution:
\begin{equation}
\label{mixedpoisson}
\prob({D^{\smallsup{G}}}=k)=\int_0^\infty e^{-x} \frac{x^k}{k!}\, dF_{\Lambda}(x),\quad k=0,1,2\ldots,
\end{equation}
where $F_{\Lambda}$ is the capacity distribution. Hence starting from
${D_1^{\smallsup{G}}},{D_2^{\smallsup{G}}},\ldots,{D_\sN^{\smallsup{G}}}$ as an i.i.d. sequence with common distribution
given by (\ref{mixedpoisson}), the (conditioned) CM with these degrees
is close to the GRG, at least in an asymptotic sense, so that one expects that
the asymptotic fluctuations of the graph distance of the CM also hold for
the generalized random graph. Also note from the mixed Poisson distribution
(\ref{mixedpoisson}), that
$$
{\tilde \nu}=\frac{\expec[D^{\smallsup{C}}(D^{\smallsup{C}}-1)]}{\expec[D^{\smallsup{C}}]}
=
\frac{\mexpec{\mL\nill^2}}{\mexpec{\mL\nill}},
$$
which is equal to $\nu$, according to \eqref{pick idd}.
As said earlier, a proof of this intuitive reasoning is missing, and our method of proof is
by coupling each random graph satisfying {\bf A1} and {\bf A2} to the Poisson random graph (PRG), and by giving a
separate proof of Theorem 1.1-1.3 for the PRG.

We finish this section by giving an overview of different distance results in random graphs.
Let $\tau$ denote the exponent of the probability mass function of the degree distribution.
In this paper and in \cite{CLaverage,finstub} the case $\tau>3$ is studied.
Results for $2<\tau<3$ for various models appeared in \cite{CLaverage,infvar,norros1,norros3}.
Typically in that case, the distance fluctuates around a constant times
$2\log\log N/|\log(\tau-2)|$.
For $1<\tau<2$,
there exists a subset of nodes with a high degree, called the \emph{core} (see \cite{infmean}).
The core forms a complete graph and almost every node is attached to the core and, thus, the graph distance is \whp\ at most $3$.

\subsection{Organization of the paper}
The  coupling argument that ties the fluctuations of the graph distance $H_\sN^{\smallsup{P}}$ in the PRG
to the fluctuations of the graph distance in random graphs satisfying assumptions {\bf A1} and {\bf A2} is treated in Section \ref{coup2}.
In Section 4 we show that the fluctuations of the graph distance $H_\sN^{\smallsup{P}}$ is given by Theorem \ref{main}.
The derivation of the fluctuations of the graph distance $H_\sN^{\smallsup{P}}$ is similar to the derivation of the fluctuations of
the graph distance $H_\sN^{\smallsup{C}}$ in the configuration model, see \cite{finstub}.
The proof in \cite{finstub} is more complicated than the proof presented here for the PRG model,
mainly because in the latter the expansion
of a given node e.g., the nodes on a given distance, can be described by means of the
so-called Reittu-Norros process, a marked branching process. This branching process will be
introduced in Section \ref{NRmodel}.

\switchText
    {Most of the proofs of the auxiliary propositions and lemmas introduced in Sections \ref{NRmodel} and \ref{sec:proof:main}
are left out as they are similar to those in \cite{finstub}. Since details of these proofs are different, we included full proofs in Appendix D of the extended version of this paper (\cite{EV}).}
    {In this paper full proofs
of the auxiliary propositions and lemmas introduced in Sections \ref{NRmodel} and \ref{sec:proof:main}
are presented in Appendix \ref{appendix coupling}. These proofs were omitted in \cite{EVshort}. }


\section{Coupling} \label{coup2}
In this section we denote by $G_\sN$ the PRG and by $G_\sN '$ some other random graph satisfying the assumptions {\bf A1} and {\bf A2},
given in Section \ref{sec:model definition}.
We number the nodes of both $G_\sN$ and $G_\sN'$ from $1$ to $N$ and we assign the  capacity $\ml i$, for $1 \leq i \leq N$, to node $i$
in each graph.
We denote by $H_{\sN}$ and $H_{\sN}'$ the graph distance between two randomly chosen nodes $A_1$ and $A_2$, such that $A_1\not = A_2$, in $G_\sN$ and $G_\sN'$, respectively.
We will show that for $N\rightarrow \infty$,
    \eq \label{Hsn=H'sn}
    \mprob{H_{\sN}\neq H_{\sN}'}=o(1).
    \en
The above implies that \whp\ the coupling of the graph distances is successful. Therefore, given the succesful coupling \eqref{Hsn=H'sn}, it is sufficient to show Theorem \ref{main} for the PRG.

The coupling bound in \eqref{Hsn=H'sn} has since been considerably
strengthened by Janson in \cite{Jans08a}, \ch{who} studies when
two random graphs with independent edges are asymptotically
equivalent. By \cite[Example 3.6]{Jans08a} and under the assumptions
{\bf A1} and {\bf A2}, we have that the
probability of {\it any} event $A_{\sN}$ is asymptotically equal for
$G_{\sN}$ and $G_{\sN}'$ \ch{when $N\rightarrow \infty$.}

\subsection{Coupling of $G_\sN$ and $G_\sN'$}   \label{subseq:couplingGG}
We next describe the coupling of the connection variables of the graphs $G_\sN$ and $G_\sN'$. A classical coupling is used, see e.g. \cite{thorisson}.
Denote by $\{X_{ij}\}_{1\leq i < j \leq N}$ and $\{X_{ij}'\}_{1\leq i < j \leq N}$
the connection variables of  the graphs $G_\sN$ and $G_\sN'$, and, similarly,
denote the connection probabilities by $\{p_{ij}\}_{1\leq i < j \leq N}$ and $\{p_{ij}'\}_{1\leq i < j \leq N}$.
For the coupling we introduce independent random variables $\{K_{ij}\}_{1\leq i < j \leq N}$.
Set $\underline p_{ij}=\min\{p_{ij},p_{ij}'\}$ and $\overline p_{ij}=\max\{p_{ij},p_{ij}'\}$, and define random variables $\hat X_{ij}$ and $\hat X'_{ij}$ with
        \begin{align*}
        \mprob{\hat X_{ij}= 1, \hat X'_{ij}= 1, K_{ij}=0}&=\underline p_{ij},&
        \mprob{\hat X_{ij}= 1, \hat X'_{ij}=0,K_{ij}=1}&=p_{ij} - \underline p_{ij},&
        \\
        \mprob{\hat X_{ij}= 0, \hat X'_{ij}=1,K_{ij}=1}&=\overline p_{ij} - p_{ij},&
         \mprob{\hat X_{ij}= 0, \hat X'_{ij}=0,K_{ij}=0}&=1-\overline p_{ij},&
        \end{align*}
        whereas all other combinations have probability 0.
        Then the laws of $\hat X_{ij}$ and $\hat X_{ij}'$ are the same as the laws of $X_{ij}$ and $X_{ij}'$, respectively. Furthermore, $K_{ij}$ assumes the value $1$ with probability $|p_{ij}-p_{ij}'|$, and is $0$ otherwise.
        Note that we do abuse the notation in the above display. We should replace the probability measure $\mathbb P$ in the above display by some other probability measure $\mathbb Q$, because the probability space is defined by the graphs $G_\sN$ and $G_\sN'$, instead of only the graph $G_\sN$. Since the graphs, conditioned on the capacities, are constructed independently from each other, this abuse of notation is not a problem.

        Consider the nodes $i$ and $j$, $1 \leq i < j \leq N$, in the graphs $G_\sN$ and $G'_\sN$ simultaneously.
        Then the event $\{K_{ij}=0\}=\{\hat X_{ij}=\hat X'_{ij}\}$ corresponds to the event that in both graphs
        there exists a connection between nodes $i$ and $j$, or that in both graphs there is no connection between nodes $i$ and $j$.
        The event $\{K_{ij}=1\}=\{\hat X_{ij} \not =\hat X'_{ij}\}$ corresponds with the event that there exists a
        connection in one of the graphs, but not in the other one.
        We call the event $\{K_{ij}=1\}$ a \emph{mismatch} between the nodes $i$ and $j$.

        Assumption {\bf A2} implies that for some constant $C' > 0$,
        \begin{align} \label{kans1}
        \mprob{K_{ij}=1}=\brc||{p_{ij}-p_{ij}'} \leq  \brc||{p_{ij}-\ml i \ml j/ l_\sN} + \brc||{p_{ij}'-\ml i \ml j/ l_\sN}
        \leq C'\frac{\ml i ^2\ml j^2}{l_\sN^2},
        \end{align}
         for all $1 \leq i < j \leq N$.
        The number of mismatches due to all the nodes incident to node $i$, $1 \leq i \leq N$, is given by
        \begin{align}\label{def:F}
        K_i=\sum_{j\not=i} K_{ij}.
        \end{align}

    Obviously, we cannot couple all the connections in the graphs $G_\sN$ and $G_\sN'$ successfully, but the total number of mismatches due to all the nodes can be bounded from above by any positive power of $N$.
To this end, we define the event ${\cal A}_\sN$ as
    \begin{align} \label{def:An}
    {\cal A}_\sN=\bigcap _{i=1}^N \bra{K_i\indic{\ml i > c_\sN}=0}= \bra{\sum_{i=1}^N K_i\indic{\ml i > c_\sN}=0},
    \end{align}
    where $c_\sN=N^{\xi}$, with $\xi > 0$.
    Then, on the event ${\cal A}_\sN$, all nodes with capacity greater than $c_\sN$ are successfully coupled.

    \begin{lemma}  \label{lemma:Ac} For each $\xi>0$ there exists a constant $\theta >0 $ such that
\begin{align} \label{eq:Ac}
    \mprob{{\cal A}_\sN^c} =\mO{ N^{-\theta}}.
\end{align}
    \end{lemma}
    \begin{proof}
We bound $\mprob{{\cal A}_\sN^c}$ using  Boole's inequality, and the
Markov inequality:
    \eqn{\label{condition:NN}
    \mprob{{\cal A}_\sN^c}
    \leq \sum_{i=1}^N \mprob{K_i\indic{\ml i > c_\sN}\geq 1}
    \leq \sum_{i=1}^N \mexpec{K_i}\indic{\ml i > c_\sN} .
    }
Then, using \eqref{kans1}, \eqref{muNnuN-def} and \eqref{reqA}, the expectation $\mexpec{K_i}$  can be bounded by,
\begin{align}
\mexpec{K_i}&=\sum_{j \neq i}
\mexpec{K_{ij}}=\sum_{j \neq i}
\mprob{K_{ij}=1}
\leq
\frac{C' \ml i^2}{l_\sN^2}\sum_{j=1}^N \ml j^2
=
\frac{C' \ml i^2}{N}\frac{\nu_\sN}{\mu_\sN}
=\mO{\lambda_i^2 N^{-1}}
\label{expec-kn}
\end{align}
Using the above and \eqref{condition:NN}, we have that
\eq
\label{condbnd}
\mprob{{\cal A}_\sN^c}
\le
\sum_{i=1}^N \mexpec{K_i}\indic{\ml i > c_\sN}
=\mO{\frac{1}{N}\sum_{i=1}^N \lambda_i^2\indic{\ml i > c_\sN}}.
\en
Observe that,
$$
\frac1N\sum_{i=1}^N \lambda_i^2\indic{\ml i > c_\sN}
=\frac1N\sum_{i=1}^N \lambda_i^{\tau-1-\varepsilon}\lambda_i^{-(\tau-3-\varepsilon)} \indic{\ml i > c_\sN}
\leq c_\sN^{-(\tau-3-\varepsilon)} \frac1N\sum_{i=1}^N\lambda_i^{\tau-1-\varepsilon} =\mO{c_\sN^{-(\tau-3-\varepsilon)}},
$$
where we applied \eqref{reqC} in the last step.
Pick $\theta=\xi(\tau-3-\varepsilon)$, then combining the above display and \eqref{condbnd}, that
$$
\mprob{{\cal A}_\sN^c}=\mO{c_\sN^{-(\tau-3-\varepsilon)}}=\mO{N^{-\theta}}.
$$
   \qed \end{proof}

   \bigskip

\subsection{Coupling the graph distances of $G_\sN$ and $G_\sN'$}
    In this subsection we couple the graph distance of the PRG with any random graph satisfying assumptions {\bf A1} and {\bf A2}.

    \begin{theorem} \label{kopHH'}
    Let $G_\sN$ be a PRG and let and $G_\sN'$ be a random graph satisfying assumption {\bf A1} and {\bf A2}.
    Let $H_{\sN}$ and $H_{\sN}'$ be the graph distances between two different uniformly chosen nodes $A_1$ and $A_2$ in, respectively, the graphs $G_\sN$ and $G_\sN'$. Then
    \eq
    \mprob{H_{\sN}\neq H_{\sN}'}=o(1).
    \en
    \end{theorem}

In order to prove Theorem \ref{kopHH'}, we will use the following strategy.
We know that for the PRG the random variable $H_\sN$ is concentrated around $\log_{\nu}{N}$. Hence if we take
$t_{\sN} =\bfloor{ (\frac 12+\eta)\log_{\nu}{N}}$, with $\eta>0$, then $\mprob{H_{\sN}>2t_{\sN}}$ is small and in order that
$H_{\sN}\neq H'_{\sN}$ on the set where $\{H_{\sN}\le 2t_{\sN}\}$, there
must have been at least one mismatch between two nodes $i$ and $j$, both being not on
a greater graph distance from our random node $A$ than  $2t_{\sN}$.

We define the neighborhood shells of a uniformly chosen node $A \in\{1,2,\ldots,N\}$, i.e., all nodes on a
fixed graph distance of node $A$, by
\begin{align} \label{def:NeiCirc}
\partial{\cal N}_0&=\{A\}&
&\text{and}&
\partial{\cal N}_k&=\{1 \leq j \leq N: d(A,j)=k\},&
\end{align}
where $d(i,j)$ denotes the graph distance between nodes $i$ and
$j$, i.e., the minimum number of edges in a path between the nodes
$i$ and $j$.
Furthermore, define the set of nodes reachable in at most $t$ steps from root $A$ by
\begin{align} \label{def:Nei}
{\cal N}_t=\{1 \leq j \leq N: d(A,j) \leq t\}=\bigcup_{k =0}^t \partial{\cal N}_k.
\end{align}
\begin{prop}
\label{prop-couplinghopc}
For $N$ sufficiently large, $t \in \mathbb{N}$, and every $b\in (0,1)$,
    \eq \label{construct}
    \mprob{H_{\sN}\neq H_{\sN}'}\leq
    \mprob{{\cal A}^c_{\sN}}+\mprob{H_{\sN}>2t}+2t\mprob{|\SPT_{t-1}|>N^b}
    +\mO{ t N^{-1+b} c_\sN^4}.
    \en
\end{prop}

\bigskip

Before giving a proof, we show that Theorem \ref{kopHH'} is a consequence of Proposition \ref{prop-couplinghopc}.

\bigskip

\noindent{\bf Proof of Theorem \ref{kopHH'}.}
By Lemma \ref{lemma:Ac}, we have that, $\mprob{{\cal A}^c_{\sN}}\leq N^{-\theta}$.
From Corollary \ref{tightness}, applied to the PRG model, we obtain that $\mprob{H_{\sN}>2t_{\sN}}=o(1)$.
The third term on the right side of \eqref{construct} can be bounded using the following lemma:
\begin{lemma} \label{LA.2.4iii}
Let $\{\SPT_t\}_{t\geq 0}$ be the reachable sets of a uniformly chosen node $A$ in the PRG $G_\sN$. Then for $\eta,\delta\in(-1/2,1/2)$ and all $t \leq (1/2+\eta)\log_\nu N$, there exists a constant $\beta_1>0$ such that
\begin{align} \label{claim.2.4}
\mprob{|\SPT_t|> N^{1/2+\delta}}= \mO{(\log_\nu N)N^{-(\delta-\eta)}}.
\end{align}
\end{lemma}
\begin{proof}
\switchText
    {The full proof of this lemma is given in the extended version of this paper \cite[Lemma D.2]{EV}.
    We now  give a heuristic derivation. We will couple $|\partial \SPT_t|$ to $\msZ t$ so that $|\partial \SPT_t|\approx \msZ t$, \whp.  If we replace $|\SPT_t|$ by $\sum_{k=0}^t \msZ k$ in \eqref{claim.2.4}, then it is easy to verify that, by using the Markov Inequality,
    $$
    \mprob{\sum_{k=1}^t \msZ k> N^{1/2+\delta}}\leq
        N^{-1/2-\delta}\sum_{k=1}^t \mexpec{\msZ k}=N^{-1/2-\delta} \sum_{k=1}^t\mu \nu^{k-1}  =\mO{ (\log_\nu N)N^{-(\delta-\eta)}},
    $$
    which then implies the result.}
    {See the proof of Lemma \ref{LA.2.4}.}
    \qed
\end{proof}

\bigskip
We now prove that all terms in the right hand of \eqref{construct} are $o(1)$ for an appropriate choice of $b$.
Lemma \ref{LA.2.4iii} implies that $2t_\sN\mprob{|\SPT_t|>N^b}=o(1)$ for some appropriately chosen $b>\frac 12$. Then, provided that $b<1$,
we see that $t_\sN N^{b-1} c_\sN^4=t_\sN N^{4\xi + b-1}=o(1)$, where we substitute $c_\sN=N^\xi$, and picking
$\xi \in (0, (1-b)/4)$.
Hence, by Proposition \ref{prop-couplinghopc}, $\mprob{H_{\sN}\neq H_{\sN}'}=o(1)$, which is precisely the content of Theorem \ref{kopHH'}.\qed

\bigskip

\noindent{\bf Proof of Proposition \ref{prop-couplinghopc}.}
We use that
    \eq
    \label{cpbd1}
    \mprob{H_{\sN}\neq H_{\sN}'}\leq
    \mprob{{\cal A}^c_{\sN}}+\mprob{H_{\sN}>2t}+
    \mprob{\{H_{\sN}\leq 2t\}\cap {\cal A}_{\sN}\cap \{H_{\sN}\neq H_{\sN}'\}}.
    \en
Let $\SPT_t^{\smallsup{i}}$ and $\SPT_t'^{\smallsup{i}}$, for $i=1,2$, be the union of neighborhood shells of the nodes $A_i$ in $G_\sN$ and $G_\sN'$, respectively.
Now, we use the fact that if $H_{\sN}\leq 2t$ and if $H_{\sN}\neq H_{\sN}'$,
then $\SPT_t^{\smallsup{1}}\not = \SPT_t'^{\smallsup{1}}$ and/or $\SPT_t^{\smallsup{2}}\not = \SPT_t'^{\smallsup{2}}$. By the exchangeability of the nodes, we have
    \eq
    \mprob{\{H_{\sN}\leq 2t\}\cap {\cal A}_{\sN}\cap \{H_{\sN}\neq H_{\sN}'\}}
    \leq 2\mprob{\{\SPT_t\neq \SPT_t'\}\cap {\cal A}_{\sN}}.
    \en
If $\SPT_t\not = \SPT_t'$, then there must be a $k\in \{1, \ldots, t\}$
for which $\SPT_k\neq \SPT_k'$, but $\SPT_{k-1}=\SPT_{k-1}'$. Thus,
    \eq \label{sumel}
    \mprob{\{H_{\sN}\leq 2t\}\cap {\cal A}_{\sN}\cap \{H_{\sN}\neq H_{\sN}'\}}
    \leq 2\sum_{k=1}^t \mprob{\{\SPT_k\neq \SPT_k'\}\cap \{\SPT_{k-1}=\SPT_{k-1}'\}\cap {\cal A}_{\sN}}.
    \en
In turn, the event $\{\SPT_k\neq \SPT_k'\}\cap \{\SPT_{k-1}=\SPT_{k-1}'\}$ implies
that one of the edges from $\partial\SPT_{k-1}$ must be miscoupled, thus $K_{i j}=1$ for some $i\in\partial\SPT_{k-1}$ and $j\in\SPT_{k-1}^c$, where $\SPT_{k-1}^c=\{1,2,\ldots,N\}\backslash \SPT_{k-1}$.
The event ${\cal A}_{\sN}$ implies that $\ml i, \ml j\leq c_\sN$.
Therefore, we bound
    \begin{multline}
    \mprob{\{\SPT_k\neq \SPT_k'\}\cap \{\SPT_{k-1}=\SPT_{k-1}'\}\cap {\cal A}_{\sN}}
    \leq \mprob{|\SPT_{k-1}|> N^b}\\
    +\sum_{i,j}\mprob{\{i\in \partial\SPT_{k-1},j\in \SPT^{c}_{k-1}, K_{ij}=1\}
    \cap \{|\SPT_{k-1}|\leq N^b\}}\indic{\ml i,\ml j\leq c_\sN}.\label{eq230}
    \end{multline}

Since $i\in\SPT_{k-1}^c$ and $j\in\partial\SPT_{k-1}$, the event $\{K_{ij}=1\}$ is
independent of $\SPT_{k-1}$ and, therefore, from $\partial\SPT_{k-1}$ as $\partial\SPT_{k-1} \subset \SPT_{k-1}$.
The edge between the nodes $i$ and $j$ points out of
$\SPT_{k-1}$, while $\SPT_{k-1}$ is determined by the occupation status of edges that are between nodes in $\SPT_{k-2}$ or pointing out of $\partial\SPT_{k-2}$. Thus, we can replace each term in the sum of \eqref{eq230} by
    \eq
    \mprob{K_{ij}=1}\mprob{\{i\in \partial\SPT_{k-1},j\in \SPT^{c}_{k-1}\}\cap\{|\SPT_{k-1}|\leq N^b\}}\indic{\ml i,\ml j\leq c_\sN}.
    \en
Since by \eqref{kans1}, we have
    \eq
    \mprob{K_{ij}=1}\indic{\ml i,\ml j\leq c_\sN}
        \leq C'
    \frac{\ml i^2\ml j^2}{l_{\sN}^2}\indic{\ml  i, \ml j\leq c_\sN}=  \mO{c_\sN^4 N^{-2}},
    \en
we can  bound the right side of \eqref{eq230} from above by
\begin{eqnarray} \nn
     \mprob{|\SPT_{k-1}|> N^b}\nn
    + \mO{c_\sN^4 N^{-2}} \sum_{i,j}\mprob{\{i\in \partial\SPT_{k-1},j\in \SPT^{c}_{k-1}\}
    \cap \{|\SPT_{k-1}|\leq N^b\}}.
    \end{eqnarray}
Finally, we  bound the sum  on the right side by
$$
 \sum_{i,j} \mprob{\{i\in \partial\SPT_{k-1},j\in \SPT^{c}_{k-1}\}\cap \{|\SPT_{k-1}|\leq N^b\}}
\leq N \mexpec{|\partial\SPT_{k-1}|\indic{|\SPT_{k-1}|\leq N^b}}
\leq N^{1+b}.
$$
Therefore, we can bound each term in the sum of \eqref{sumel} by
    \begin{eqnarray} \nn
    \shift &&\mprob{\{\SPT_k\neq \SPT_k'\}\cap \{\SPT_{k-1}=\SPT_{k-1}'\}\cap {\cal A}_{\sN}}
    \leq \mprob{|\SPT_{k-1}|> N^b} + \mO{c_\sN^4 N^{-1+b}}.
    \end{eqnarray}
Since, for $k \leq t$, we have that $\mprob{|\SPT_{k-1}| > N^b } \leq \mprob{|\SPT_{t-1}| > N^b }$, by summing over $k=1,\dots, t$ in \eqref{sumel},
 we arrive at
   $$
    \mprob{\{H_{\sN}\leq 2t\}\cap {\cal A}_{\sN}\cap \{H_{\sN}\neq H_{\sN}'\}}
    \leq 2 t \mprob{|\SPT_{t-1}| > N^b }
    +\mO{t c_\sN^4N^{-(1-b)}}.
    $$
Together with \eqref{cpbd1} this proves the proposition.
\qed

\section{The Poissonian random graph model} \label{NRmodel}
The proof of the fluctuations of the graph distance in the CM in \cite{finstub} has been done in a number of steps.
One of the most important steps is the coupling of the expansion of the neighborhood shells  of a node to  a BP.
For the PRG, we follow the same strategy as in \cite{finstub}, although the details differ substantially.

The first step is to introduce the \mRNprocess, which is a marked BP.
The \mRNprocess\ was introduced by Norros and Reittu in \cite{norros3}.
We can thin the \mRNprocess\ in such a way that the resulting process, the \mrRNprocess,
can be coupled to the expansion of the neighborhood shells of a randomly chosen node in the PRG.
Finally, we introduce capacities for the \mRNprocess\ and the \mrRNprocess.

\bigskip
\subsection{The NR-process and its thinned version} \label{subsec:static} The \mRNprocess\ is
a marked delayed BP denoted by $\{\mZ t,\mJ
t\}_{t\ge 0}$, where $\mZ t$ denotes the number of individuals of
generation $t$, and where the vector
$$
\mJ t=(\mM[t] 1,\mM[t] 2,\ldots,\mM[t] {\mZ t}) \in \{1,2,\ldots,N\}^{\mZ t},
$$
denotes the marks of the individuals in generation $t$.
We now give a more precise definition of the \mRNprocess\ and describe its connection with  the PRG.
We define $\mZ 0=1$ and take $\mM[0] 1$ randomly from the set $\{1,2,\ldots,N\}$,
corresponding to the choice of $A_1$.
The offspring of an individual with mark $m\in \{1,2,\ldots,N\}$ is as
follows: the total number of children has a Poisson distribution
with parameter $\ml m$, of which, for each $i\in \{1,2,\ldots,N\}$, a Poisson
distributed number with parameter
\begin{equation}
\label{parameter} \frac{\ml i\ml m}{l_\sN},
\end{equation}
bears mark $i$, independently of the other individuals. Since
$\sum_{i=1}^N \ml i\ml m/l_\sN=\ml m$,
and sums of independent Poisson random variables are again
Poissonian, we may take the number of children with different marks mutually independent. As a result of this definition, the
marks of the children of an individual in $\{\mZ t,\mJ t\}_{t\ge 0}$ can be seen as independent realizations of a random variable
$M$, with distribution
\begin{equation}
\label{dist mark} \mprob{M=m}= \frac{\ml m} {l_\sN}, \quad
 1 \leq m \leq N,
\end{equation}
and, consequently,
\begin{align}\label{ml M}
\mexpec{\ml M}=
\sum_{m=1}^N\ml m\mprob{M=m} =
\frac{1}{l_\sN}\sum_{m=1}^N\ml m^2.
\end{align}

For the definition of the \mrRNprocess\ we start with a copy of the \mRNprocess\
$\{\mZ t, \mJ t\}_{t\geq0}$, and
reduce this process generation by generation, i.e., in the order
\begin{align} \label{full seq}
\mM[0] 1,\mM [1] 1 ,\ldots \mM [1]{\mZ 1},\mM [2]{1},\ldots
\end{align}
by discarding each individual and all its descendants whose mark
has appeared before.
The process obtained in this way is called the  \mrRNprocess\ and
is denoted by the sequence $\{\mnZ t,\mnJ t\}_{t\geq0}$.
One of the main results of \cite{norros3} is Proposition 3.1:

\begin{prop} \label{N NR kop}
Let $\{\mnZ t,\mnJ t\}_{t\geq 0}$ be the \mrRNprocess\ and let
$\mnJ t$ be the set of marks in the $t{\rm -th}$ generation, then the sequence of sets
$\{\mnJ t\}_{t\ge 0}$
has the same distribution as the sequence $\{\partial{\cal
N}_t\}_{t\geq0}$ given by \eqref{def:NeiCirc}.
\end{prop}
\begin{proof}
In \cite{norros3} it is assumed that the sequence of capacities is random, but in most proofs, including \cite[Proposition 3.1]{norros3}, Norros and Reittu condition on the capacities, and therefore consider the capacities as deterministic. Thus, the proof of \cite[Proposition 3.1]{norros3} holds verbatim. \qed
\end{proof}

\bigskip

As a consequence of the previous proposition, we can couple the \mrRNprocess\ to the
neighborhood shells of a uniformly chosen node $A\in\{1,2,\ldots,N\}$, i.e., all nodes on a
fixed graph distance of $A$, see \eqref{def:NeiCirc} and note that $A \sim \mM[0] 1$.
Thus, using the the above proposition, we can couple the expansion of the neighborhood shells and the \mrRNprocess\ in such a way that
\begin{align} \label{kop nr ne}
\mnJ t&=\partial{\cal N}_t& &\text{and}& \mnZ t &= |\partial{\cal N}_t|, \, \, t\geq 0.&
\end{align}
Furthermore, we see that an individual with mark $m$, $1\leq m\leq N$, in the \mrRNprocess\ is identified with node $m$ in the graph $G_\sN$,
whose capacity is given by $\ml m$.

The offspring distribution $\mf\nill$ of $\mZ 1$, i.e., the first generation of $\{\mZ t\}_{t\geq 0}$, is given by
\begin{align} \label{def mf k}
\mf n  &=
 \mprob{\mP[\ml{A}] = n}
 =\frac{1}{N}\sum_{m=1}^N e^{-\ml m}\frac{\ml m^n}{n!},\qquad n \geq 0.
\end{align}
Recall that individuals in the second and further generations have a random mark distributed as  $M$,
given by \eqref{dist mark}.
Hence, if we denote the offspring distribution of the second and further generations by $\mg n$, then we obtain
\begin{align} 
\mg n &=\mprob{{\mP[\ml M]=n}}
= \sum_{m=1}^N e^{-\mL m} \frac{\ml m^{n}}{n!}\frac{\ml m}{l_\sN}
\label{bepaling g}
    = \frac{1}{l_\sN}\sum_{m=1}^N e^{-\ml m} \frac{\ml m^{n+1}}{n!}, \qquad n\ge 0.
\end{align}
Furthermore, we can relate $\mg n$ and $\mf n$ by
\begin{align} \label{relatie f en g}
\mg n
=\frac{(n+1)}{l_\sN/N}\frac{1}{N}\sum_{m=1}^N e^{-\ml{m}}
\frac{\ml{m}^{n+1}}{(n+1)!} =\frac{(n+1)\mf {n+1}}{\mu_\sN}.
\end{align}
It follows from condition (\mCo2) that $\mf n\to f_n$ and $\mg n\to g_n$.



\subsection{Coupling  with a delayed BP} \label{subsec:BPcoupling}

In this subsection we  will introduce a coupling between the \mrRNprocess\ and the delayed BP $\{\msZ t\}_{t\geq0}$, which is defined by condition (\mCo2) in Section \ref{sec: main results}.
This coupling is used in the proof of Theorem \ref{main} and \ref{limit law} for the PRG, to express the probability distribution of $H_\sN$ in terms of the BP $\{\msZ t\}_{t\geq0}$.
\switchText{The full proof of these propositions are given in the extended version of this paper \cite{EV}.}{}

Introduce the total capacity of the $t$-th generation of the \mRNprocess\ $\{\mZ t, \mJ t\}_{t\geq0}$ and the \mrRNprocess\ $\{\mnZ t, \mnJ t\}_{t\geq0}$ as, respectively,
\begin{align} \label{def C&Ct}
\mC {t+1}&=\sum_{i=1}^{\mZ t} \mli{\mM[t] i}&
&\text{and}&
\mnC {t+1}&=\sum_{i=1}^{\mnZ t} \mli{\mnM[t] i}, \,\, t\geq 0,&&
\end{align}
where, to improve readability, we write $\lambda(A)=\lambda_A$.
Using the coupling given by \eqref{kop nr ne}, we can rewrite the capacity $\mnC {t+1}$ as
\begin{align}\label{def CinN}
\mnC {t+1}=\sum_{i\in\partial\SPT_t} \ml i.
\end{align}
For the proof of Theorem \ref{main} and \ref{limit law}, in the case of the PRG, we need to control the difference between $\mC t$  and $\mnC t$ for fixed $t$.
For this we will use the following proposition:
\begin{prop} \label{coupling of NR}
There exist constants $\alpha_2,\beta_2>0$, such that for all $0<\eta <\alpha_2$ and
all $t\leq t_{\sN}=(1/2+\eta)\log_\nu N$,
\begin{align}
\mprob{\sum_{k=1}^t(\mC k -\mnC k) > N^{1/2-\alpha_2}} \leq N^{-\beta_2}.
\end{align}
\end{prop}
\begin{proof}
\switchText
    {
    Notice that $\mnC k \leq \mC k$ holds trivially, because $\mnZ k$ is obtained from $\mZ k$ by thinning.
    The full proof of Proposition \ref{coupling of NR} is given in Section D.1 of \cite{EV}.
    The proof consists of several steps.
    Denote by a \emph{duplicate} an individual in the \mRNprocess\ whose mark has appeared previously.
    In  \cite{EV} we will give a formal definition of a duplicate.
    Firstly, we have to keep track of all the duplicates.
    We will show that \whp\ duplicates do not appear in the first $t_{\sN}$ generations and that the number of duplicates in the first $t_{\sN}$ generations can be bounded from above by some small power of $N$. Then, secondly, we will bound the total progeny of each duplicate. Combining these results  gives the claim of this proposition.}
    {The proof is deferred to Section \ref{coupling of NR:A}.}
\qed
\end{proof}

\bigskip

In order to prove Theorem \ref{main} and Theorem \ref{limit law} we will grow two \mrRNprocess es $\{\mnZi i t, \mnJi i t\}_{t\geq 0}$, for $i=1,2$.
The root  of $\{\mnZi i t, \mnJi i t\}_{t\geq 0}$ starts from a uniformly chosen node or mark $A_i \in \bra{1,2,\ldots,N}$. These two nodes are different \whp, because
$$
\mprob{A_1=A_2}=\frac{1}{N}.
$$
By \eqref{kop nr ne} the \mrRNprocess\ can be coupled to the neighborhood expansion shells $\{\SPT_t^{\smallsup{1}}\}_{t\geq 0}$ and  $\{\SPT_t^{\smallsup{2}}\}_{t\geq 0}$. In the following lemma we compute the distribution of the number of edges between two shells with different subindeces, i.e., $\SPT_k^{\smallsup{1}}$ and $\SPT_t^{\smallsup{2}}$.

\begin{lemma} \label{lemma:distribution of edges}
Fix integers $k$ and $t$. Then conditionally on $\SPT_k^{\smallsup{1}}$ and $\SPT_t^{\smallsup{2}}$ and
given that $\SPT_k^{\smallsup{1}} \cap \SPT_t^{\smallsup{2}} = \emptyset$ the number of edges between the nodes in
$\SPT_k^{\smallsup{1}}$ and $\SPT_t^{\smallsup{2}}$
is distributed as a Poisson random variable with mean
\begin{align}
\frac{\mnCi 1 {k+1} \mnCi 2 {t+1}}{l_\sN}. \label{mcimcj}
\end{align}
\end{lemma}
\begin{proof}
Conditioned on $\SPT_k^{\smallsup{1}}$, $\SPT_t^{\smallsup{2}}$ and $\SPT_k^{\smallsup{1}} \cap \SPT_t^{\smallsup{2}} = \emptyset$, the number of edges between $\SPT_k^{\smallsup{1}}$ and $\SPT_t^{\smallsup{2}}$ is given by
\begin{align} \label{poissonSum}
\sum_{i\in \partial\SPT_k^{\smallsup{1}}}\sum_{j\in \partial\SPT_t^{\smallsup{2}}} E^{\smallsup{P}}_{ij},
\end{align}
where $E_{ij}^{\smallsup{P}}$ are independent Poisson random variables with mean $\ml i\ml j/l_\sN$, see \eqref{poissonedge}. Sums of independent Poisson random variables are again Poissonian, thus \eqref{poissonSum} is a Poisson random variable with mean the expected value of \eqref{poissonSum}:
$$
\sum_{i\in \partial\SPT_k^{\smallsup{1}}}\sum_{j\in \partial\SPT_t^{\smallsup{2}}} \frac{\ml i\ml j}{l_\sN}
=\frac{\mnCi 1 {k+1} \mnCi 2 {t+1}}{l_\sN},
$$
where we have used \eqref{def CinN} in the last step. \qed
\end{proof}
\bigskip

The further proof of Theorems \ref{main}-\ref{limit law} crucially relies on the following technical claim:
\begin{prop} \label{coupling of sums}
There exist constants $u_2,v_2,\eta>0$ such that for all $t\le t_{\sN}=(1+2\eta)\log_{\nu} N$,
as $N\to\infty$,
\begin{eqnarray} \label{coupling of sums1:eq}
\prob\Big(\frac{1}{N}\Big|
\sum_{k=2}^{t+1} \msZi 1 {\lceil k/2 \rceil} \msZi 2 {\lfloor k/2 \rfloor}
-
\sum_{k=2}^{t+1}\mnCi 1 {\lceil k/2 \rceil} \mnCi 2 {\lfloor k/2\rfloor}\Big|>N^{-u_2}\Big)=\mO{N^{-v_2}}.
\end{eqnarray}
\end{prop}
\begin{proof}
\switchText
    {The proof of Proposition \ref{coupling of sums} is given in the extended version of this paper, see \cite{EV}.
    We will intuitively explain how $\mnCi i k$, for $i=1,2$ and $1\leq k \leq t$, can be replaced by $\msZi i k$. Firstly, whenever $k \leq (1/2+\eta)\log_\nu N$, we can neglect the influence of the thinning.
    Thus, we can replace $\mnCi i k$ by $\mCi i k$.
  The capacity $\mCi i k$ is the sum of the capacities of the $\mZi i {k-1}$ nodes of the $(k-1)$-th generation.
Conditionally on $\mZi i {k-1}$, the value of $\mexpec{\mCi i k|\mZi i {k-1}}$ is given by $\nu \mZi i {k-1}$, thus the next step is to replace $\mCi i k$ by $\nu \mZi i {k-1}$.
    Finally, by condition (\mCo2) the offspring distributions of the BP $\{\mZ k\}_{k\geq 0}$ converges to those of $\{\msZ k\}_{k\geq 0}$. Therefore, we can replace $\nu \mZi i {k-1}$ by $\nu \msZi i {k-1}$, which we, finally, replace by $\msZi i {k}$.}
    {The proof is deferred to Section \ref{proposition:coupling of sums:bewijs}.}
\qed
\end{proof}

\bigskip

In the next section we will use this proposition in combination with Lemma \ref{lemma:distribution of edges} to replace sums over capacities,
which do not depend on $N$, by sums over sizes of a BP, which no longer depend on $N$.


\section{Proof of Theorem \ref{main} and \ref{limit law} for the PRG} \label{sec:proof:main}
In this section, we prove Theorem \ref{main} and \ref{limit law} for the PRG model.
Using the coupling result in Proposition \ref{prop-couplinghopc} we obtain Theorem \ref{main} and \ref{limit law} for all random graphs satisfying the assumptions {\bf A1} and {\bf A2}.
As in the previous section, we denote by $G_\sN$ a PRG.

We grow two \mrRNprocess es. Each \mrRNprocess\ starts from a uniformly chosen node $A_i\in\{1,2,\ldots,N\}$, $i=1,2$, such that $A_1\not=A_2$, \whp.

\bigskip
\noindent{\bf Step 1: Expressing $\prob(H_\sN>t)$ in capacities.}
We have $H_\sN>1$ iff (if and only if) there are no edges
between the nodes $A_1$ and $A_2$. Given the capacities
$\mnCi 1 1$ and $\mnCi 2 1$, the number of edges between the nodes $A_1$ and $A_2$
has, according to Lemma \ref{lemma:distribution of edges}, a Poisson distribution with mean
\begin{equation}
\label{poissoncc}
\frac
{\mnCi 1 1 \mnCi 2 1}
{l_\sN}.
\end{equation}
Therefore,
\begin{equation}
\label{kansH>1}
\prob(H_\sN>1)
=\mexpec{\exp\left\{-\frac
{\mnCi 1 1\mnCi 2 1}
{l_\sN}\right\}}.
\end{equation}

We next inspect the capacity of the first generation of $\mnZi 1 1$ given by $\mnCi 1 2$.
Given $\mnCi 1 2$ and $\mnCi 2 1$, that is the total capacity
of the nodes in $\mnZi 1 1$ and the capacity of node $A_2$, we again have a Poisson
number of edges between node $A_2$ and the nodes in $\mnZi 1 1$,
however, this time with parameter
\begin{equation}
\label{poissonc2}
\frac
{\mnCi 1 2 \mnCi 2 1}
{l_\sN}.
\end{equation}

In order to compute the survival probability $\prob (H_\sN>t)$ we need more notation.
We write $\Q_C^{\smallsup{t_1,t_2}}$ for the conditional probabilities given
$\{\mnCi 1 k\}_{k=1}^{t_1}$ and $\{\mnCi 2 k\}_{k=1}^{t_2}$. We further write
$\expec_C^{\smallsup{t_1,t_2}}$ for the expectation with respect to $\Q_C^{\smallsup{t_1,t_2}}$.
For $t_2=0$, we only condition on
$\{\mnCi 1 k\}_{k=1}^{t_1}$.
Lemma \ref{lemma:distribution of edges} implies that
$$
\Q_C^{\smallsup{k,t}}(H_\sN>k+t-1|H_\sN > k+ t -2)=\exp\bra{-\frac{\mnCi 1 {k} \mnCi 2 {t} }{l_\sN}}.
$$
Indeed, the event $\{H_\sN > k+ t -2\}$ implies that $\SPT_{k}^{\smallsup{1}} \cap \SPT_{t}^{\smallsup{2}} = \emptyset$.
From (\ref{poissoncc}) and the above statement,
\begin{eqnarray*}
\prob(H_\sN>2)
&=&\mexpec{\Q_C^{\smallsup{1,1}}(H_\sN>1)\,\Q_C^{\smallsup{1,1}}(H_\sN>2|H_\sN>1)}\\
&=&\mexpec{\Q_C^{\smallsup{1,1}}(H_\sN>1)\,\expec_C^{\smallsup{1,1}}
     \brh{Q_C^{\smallsup{2,1}}(H_\sN>2|H_\sN>1)}}\\
     &=&
     \mexpec{
        \expec_C^{\smallsup{1,1}}
                \brh{\Q_C^{\smallsup{1,1}}(H_\sN>1)\,\Q_C^{\smallsup{2,1}}(H_\sN>2|H_\sN>1)
                }} \\
     &=&\mexpec{\exp\left\{-\frac
{\mnCi 1 1 \mnCi 2 1}
{l_\sN}\right\}\cdot \exp\left\{-\frac
{\mnCi 1 2 \mnCi 2 1}
{l_\sN}\right\}}
\\&=&
\mexpec{
        \exp\left\{
   - \frac{\sum_{k=2}^{3} \mnCi 1 {\lceil k/2 \rceil}
    \mnCi 2 {\lfloor k/2\rfloor}} {l_\sN} \right\}
    }.
\end{eqnarray*}
By induction we obtain as in \cite[Lemma 4.1]{finstub},
\begin{align} \label{compeven}
\prob(H_\sN>t)= \mexpec{
    \exp\left\{
    -\frac{\sum_{k=2}^{t+1} \mnCi 1 {\lceil k/2 \rceil}
    \mnCi 2 {\lfloor k/2\rfloor}} {l_\sN} \right\}
    }.
\end{align}
Note that \eqref{compeven} is an equality, while in \cite{finstub} an error needed to be taken along.

\bigskip

\noindent {\bf Step 2: Coupling with the delayed BP.}
In this step we replace $\mnCi 1 t$ and $\mnCi 2 t$ by $\msZi 1 t$ and $\msZi 2 t$.

For each event $\cal B$, and any two nonnegative random variables $V$ and $W$,
\begin{align*}
\brc||{\mexpec{e^{-V}}-\mexpec{e^{-W }}}
    &\leq \brc||{\mexpec{(e^{-V}-e^{-W})\mindic{\cal B}}} + \mprob{{\cal B}^c}.
\nonumber
\end{align*}
Now take
$$
{\cal B}=\left\{\frac{1}{N}\Big|\sum_{k=2}^{t+1}
\msZi 1 {\lceil k/2 \rceil} \msZi 2 {\lfloor k/2
\rfloor}- \sum_{k=2}^{t+1}\mnCi 1 {\lceil k/2 \rceil}
\mnCi 2 {\lfloor k/2\rfloor}\Big|\leq N^{-u_2}\right\},
$$
and the random variables $V$ and $W$ as
$$
V=\frac{1}{N}
\sum_{k=2}^{t+1} \msZi 1 {\lceil k/2
\rceil} \msZi 2 {\lfloor k/2\rfloor} ,
\qquad W=
\frac{1}{N}\sum_{k=2}^{t+1} \mnCi 1 {\lceil k/2\rceil}
\mnCi 2 {\lfloor k/2\rfloor}.
$$
Then for $t \leq t_{\sN}$, the Proposition \ref{coupling of sums} implies that $\prob({\cal B}^c)=\mO{N^{-v_2}}$, whereas on the event ${\cal B}$ we have $|V-W|\le N^{-u_2}$.
Hence, using that $e^{-v}-e^{-w}=\mO{v-w}$ when $v,w$ are small, and that $e^{-v}\le 1,\, v\ge 0$,
we obtain
\eq
\label{verschil-exp}
\brc||{\mexpec{e^{-V}}-\mexpec{e^{-W }}}
\le
\mO{N^{-u_2}}\prob({\cal B})+\prob({\cal B}^c)=\mO{N^{-u_2}}+\mO{N^{-v_2}}.
\en
It is now clear from {\it step} 1, the above result and by replacing $l_\sN$ by $N\mu$, using \eqref{reqA},
that for some $\beta>0$ and all $t \leq t_{\sN}$,
\begin{align}\label{hhoednaarexp}
    \prob(H_\sN>t)=
        \mexpec{
            \exp\left\{-
            \frac{\sum_{k=2}^{t+1} \msZi 1 {\lceil k/2 \rceil} \msZi 2 {\lfloor
            k/2\rfloor} } {\mu N} \right\}
            } +\mO{N^{-\beta}}.
\end{align}

\bigskip

\noindent {\bf Step 3: Evaluation of the limit points.}
From this step on, the remainder of the proof of our main theorem is identically to the proof of Theorem 1.1 in \cite{finstub}.
To let the proof be self-contained, we finish the main line of the argument.
Starting from (\ref{hhoednaarexp}) we replace $t$ by $\sigma_{\sN}+t$ and assume that $\sigma_{\sN}+t\le
t_{\sN}$, where, as before, $\sigma_{\sN}=\lfloor\log_{\nu}
N\rfloor$, to obtain
    \begin{equation}
    \label{biggerthanlogN}
    \prob(H_\sN>\sigma_{\sN}+t)=
    \mexpec{
        \exp\left\{
        -\frac{\sum_{k=2}^{\sigma_{\sN}+t+1} \msZi 1 {\lceil k/2
        \rceil} \msZi 2 {\lfloor
        k/2\rfloor} } {\mu N} \right\}
    }   +\mO{N^{-\beta}}.
    \end{equation}
We write $N=\nu^{\log_{\nu}N}=\nu^{\sigma_{\sN}-a_{\sN}},$ where
we recall that $a_{\sN}=\lfloor \log_{\nu}N\rfloor-\log_{\nu}N$.
Then
$$
\frac{\sum_{k=2}^{\sigma_{\sN}+t+1} \msZi 1 {\lceil k/2
\rceil} \msZi 2 {\lfloor k/2\rfloor}}{\mu N} = \mu \nu^{a_{\sN}+t}
\frac{\sum_{k=2}^{\sigma_{\sN}+t+1} \msZi 1 {\lceil k/2
\rceil} \msZi 2 {\lfloor k/2\rfloor}}{\mu^2 \nu^{\sigma_{\sN}+t}}.
$$
In the above expression, the factor $\nu^{a_{\sN}}$ prevents
proper convergence.

Without the factor $\mu \nu^{a_{\sN}+t}$, we
obtain from Appendix A4 of \cite{finstub} that, with probability $1$,
    \begin{equation}
    \label{convergingemiddelde}
    \lim_{N\to \infty} \frac{\sum_{k=2}^{\sigma_{\sN}+t+1}
    \msZi 1 {\lceil k/2 \rceil} \msZi 2 {\lfloor  k/2\rfloor}}
    {\mu^2 \nu^{\sigma_{\sN}+t}} =\frac{{\cal
    W}^ {\smallsup{1}}{\cal W}^{\smallsup{2}}}{\nu-1},
    \end{equation}
We now use the speed of convergence result of \cite{asmussen}, which was further
developed in Section 2 of \cite{finstub} and which reads that there exists a positive $\beta$
such that:
\eq
\label{big-O}
\prob(|{\cal W}-{\cal W}_k|>(\log N)^{-\alpha})=\mO{N^{-\beta}}, \, k \leq \bfloor{\frac12\log_\nu N}.
\en
for each positive $\alpha$.
Combining (\ref{biggerthanlogN})  and (\ref{big-O}) we obtain that for each $\alpha>0$ and for $t\le 2\eta\log_{\nu}N$,
\eq
\label{limietresultaat}
\prob(H_\sN>\sigma_{\sN}+t)
=
\mexpec{\exp\{-\kappa\nu^{a_{\sN}+t}{\cal W}_1{\cal W}_2\}}+\mO{(\log N)^{-\alpha}}.
\en
From \eqref{limietresultaat} one can finally derive as in \cite{finstub}, that, asymptotically as
$N\to \infty$, the probability $\prob(H_\sN<\infty)$ is
equivalent to the probability $q^2=\prob({\cal W}_1{\cal W}_2>0) $, where $q$ is the survival probability of the branching process $\{\msZ t\}_{t\geq0}$, so that
(\ref{limietresultaat}) induces for $t\le 2\eta\log_{\nu}N$,
\eq
\label{condconver}
\prob(H_\sN\leq \sigma_{\sN}+t|H_\sN<\infty)
=
\mexpec{1-\exp\{-\kappa\nu^{a_{\sN}+t}{\cal W}_1{\cal W}_2\}|{\cal W}_1{\cal W}_2>0}+o(1).
\en

%
%
%
%
%
\newpage
\renewcommand{\thesection}{\Alph{section}}
\setcounter{section}{0}

\numberwithin{equation}{section}
\numberwithin{theorem}{section}

\newcommand{\mvar}[1]{\text{Var}\br{#1}}
\section{Appendix: Mixed Poisson distributions}
\ch{In this section, we study mixed Poisson distributions.}

\begin{lemma}\label{lem.mixingpoisson}
    Consider a mixed Poisson distribution with mixing random variable $\Theta_f$. Let, for $n=0,1,\ldots$,
    $$
    f_n=\mprob{\mPoisson{\Theta_f}=n},
    $$
    and consider the sequence $\{g_n\}_{n\geq 0}$, where
    $$
    g_n=\frac{(n+1)f_{n+1}}{\mexpec{\Theta_f}},
    $$
    then $\{g_n\}_{n\geq 0}$ is also a mixed Poisson distribution with mixing random variable $\Theta_g$, where
    \eqn{\label{A.mi.3}
    \mprob{\Theta_g \leq x} =\frac{1}{\mexpec{\Theta_f}}\mexpec{\Theta_f\indic{\Theta_f \leq x}}.
    }
\end{lemma}
\begin{proof}
    We will assume that \eqref{A.mi.3} holds and show that $\mprob{\mPoisson{\Theta_g}=n}=g_n$.
    Let $F_f$ and $F_g$ be the distribution functions of $\Theta_f$ and $\Theta_g$, respectively.
    Observe that
    $$
        \frac{\mathrm{d}F_g(x)}{\mathrm{d}x}
                =\frac x {\mexpec{\Theta_f}} \frac{\mathrm{d}F_{f}(x)}{\mathrm{d}x}.
    $$
    Therefore,
    \begin{align*}
    \mprob{\mPoisson{\Theta_g}=n}
    &=\int_0^\infty e^{-x}\frac{x^n}{n!}\,\mathrm dF_g(x)
    =\frac 1 {\mexpec{\Theta_f}}\int_0^\infty e^{-x}\frac{x^{n+1}}{n!}\,\mathrm dF_f(x)
    \\&=\frac {n+1} {\mexpec{\Theta_f}}\int_0^\infty e^{-x}\frac{x^{n+1}}{(n+1)!}\,\mathrm dF_f(x)=\frac{(n+1)f_{n+1}}{\mexpec{\Theta_f}}=g_n,
    \end{align*}
    which proves the claim. \qed
\end{proof}

\bigskip

As an example, we consider the random variable
$\Theta_{f}^\smallsup\sN$ that takes the value $\lambda_i$,
$i=1,2,\ldots,N$, with probability $p_i,$
then, using Lemma \ref{lem.mixingpoisson}, $\Theta_{g}^\smallsup\sN$ takes the value $\lambda_i$, $i=1,2,\ldots,N$, with probability
    \eqn{\label{mixpoissofg}
    \mprob{\Theta_g^{\smallsup{N}}=\lambda_i}=\frac{p_i\lambda_i}{\sum_{j=1}^Np_j\lambda_j}.
    }

\ch{We close this section with a lemma that relates condition (\mCo3)
to a condition on the capacities $\{\lambda_i\}_{i=1}^N$:}

\begin{lemma}[Moments of $\{\mg n\}_{n\geq0}$ and $\{\lambda_i\}_{i=1}^N$]
\label{lem-gvslambda}
Let $q\geq 1$, and assume that (\mCo1) holds. Then,
    \eqn{\label{event-gvslambda}
    \limsup_{N\rightarrow \infty}\sum_{n=0}^{\infty} n^q \mg n<\infty
    \qquad\ch{\text{if and only if}}\qquad
    \limsup_{N\rightarrow \infty}\frac 1N \sum_{i=1}^{N} \lambda_i^{q+1}<\infty.
    }
\end{lemma}

\bigskip

\proof Observe from \eqref{def-fgn-n}, that
    \eqn{
    \label{momgnrew}
    \sum_{n=0}^{\infty} n^q \mg n
    =\frac{1}{N\mu_{\sN}}\sum_{i=1}^N \lambda_i \sum_{n=0}^{\infty}  n^{q} e^{-\lambda_i}\frac{\lambda_i^{n}}{n!}
    =\frac{1}{N\mu_{\sN}}\sum_{i=1}^N \lambda_i \expec\big[\mPoisson{\lambda_i}^{q}\big].
    }
Since $q\geq 1$, we have that, by the H\"older inequality,
    \eqn{
    \expec\big[\mPoisson{\lambda_i}^{q}\big]\geq \expec\big[\mPoisson{\lambda_i}\big]^q
    =\lambda_i^q.
    }
As a result, we obtain that
    \eqn{
    \frac 1N \sum_{i=1}^N \lambda_i^{q+1}\leq \mu_{\sN} \sum_{n=0}^{\infty} n^q \mg n.
    }
By (\mCo1), $\limsup_{N\rightarrow \infty}\mu_{\sN} =\mu<\infty$, so that
indeed $\limsup_{N\rightarrow \infty}\frac 1N \sum_{i=1}^{N} \lambda_i^{q+1}<\infty$ follows
when $\limsup_{N\rightarrow \infty}\sum_{n=0}^{\infty} n^q \mg n<\infty$. This proves
the first implication.

For the second implication, we make use of the fact that there exists
a constant $c>0$ such that, uniformly for all $\lambda\geq 0$,
$\expec\big[\mPoisson{\lambda}^{q}\big]\leq c(\lambda^{q}+1).$
Then, by \eqref{momgnrew}, we have that
    \eqn{ \label{helper:mom}
    \sum_{n=0}^{\infty} n^q \mg n
    \leq \frac{c}{N\mu_{\sN}}\sum_{i=1}^N \lambda_i(\lambda_i^{q}+1)
    =c+\frac{c}{\mu_{\sN}}\frac 1N\sum_{i=1}^N \lambda_i^{q+1}.
    }
By (\mCo1), $\liminf_{N\rightarrow \infty}\mu_{\sN} =\mu>0$, so that
when $\limsup_{N\rightarrow \infty}\frac 1N \sum_{i=1}^{N} \lambda_i^{q+1}<\infty$,
we obtain
    \eqn{
    \limsup_{N\rightarrow \infty}\sum_{n=0}^{\infty} n^q \mg n
    \leq c+\frac{c}{\mu} \limsup_{N\rightarrow \infty} \frac 1N \sum_{i=1}^{N} \lambda_i^{q+1}
    <\infty.
    }
This proves the second implication, and hence completes the proof.
\qed

\section{Appendix: Deterministic capacities}
\label{app_detcap}
\newcommand{\sinv}{{\bar F}^{-1}}
\newcommand{\plaf}{\lceil NU \rceil}

\ch{Let $F$ be any distribution function, and write $\bar{F}(x)=1-F(x)$
for the survival function of $F$.
Define the inverse $\bar{F}^{-1}(u)$, for $u\in (0,1)$, by
    \begin{equation}
    \label{invverd}
    \sinv(u)=\inf \{ s: \bar{F}(s)\leq u\},
    \end{equation}
where the definition is chosen such that
    \eqn{
    \label{Finvid}
    \sinv(1-u)=F^{-1}(u)=\inf\{x: F(x)\geq u\}.
    }
We shall often make use of \eqref{Finvid}, in particular since it
implies that $\sinv(U)$ has distribution function $F$ when $U$ is
uniform on $(0,1)$.} Throughout this section we use the
abbreviations $u_{\sss N}=\lceil uN\rceil /N$ and $U_{\sss N}=\lceil UN\rceil /N$.

In this section, we let the capacities $\{\lambda_i\}_{i=1}^N$ be
given by
    \eq
    \label{deflam}
    \lambda_i=\sinv(i/N), \qquad i=1,2,\ldots,N,
    \en
where $F\colon [0,\infty)\mapsto [0,1]$ is a distribution function
satisfying \eqref{distribution}. \ch{We shall verify that conditions (\mCo1)--(\mCo3)
hold under appropriate conditions on $F$, the main condition being \eqref{distribution},
but we shall need one more technical assumption.}

We let the random variable
$\Theta^{\smallsup{N}}_f$ takes the value $\lambda_i$,
$i=1,2,\ldots,N$, with probability $p_i=\frac1N$, and note that
the distribution of $\Theta^{\smallsup{N}}_f$ equals
$\Theta^{\smallsup{N}}_f=\sinv(U_{\sss N})$, where $U$ has a uniform
$(0,1)$ distribution. We start by identifying
\ch{$f, f^{\smallsup{N}}$, $g$ and $g^{\smallsup{N}}$.}
We observe that, in this setting,
    $$
    \mf n=\expec[\prob(Poi(\Theta^{\smallsup{N}}_f=n)]=\expec\left[e^{-\sinv\Big(U_{\sss N}\Big)}
    \frac{\sinv\Big(U_{\sss N}\Big)^n}{n!}\right].
    $$
We define
    \eq
    \label{deffn}
    \mfl n
    =\expec\left[e^{-\sinv(U)}
    \frac{\sinv(U)^n}{n!}\right]
    =
    \expec[\prob(\mPoisson{\Theta_f)=n}],
    \en
where $\Theta_f=\sinv(U)$.
Observe that the used $U$ is the same in the definition of $\Theta^{\smallsup{N}}_f$, thus we implicitly introduced a coupling between
the random variables $\Theta_f$ and $\Theta^{\smallsup{N}}_f$.
Then, using that
$\Theta_f>\Theta^{\smallsup{N}}_f$ a.s.,\ \ch{since $u\mapsto \sinv(u)$ is non-increasing},
    \begin{align}\nonumber
    \frac12 \sum_{n=1}^\infty |\mf n-\mfl n|&= \mbox{d}_{\sss \rm TV}(Poi(\Theta^{\smallsup{N}}_f),Poi(\Theta_f))\leq
    \prob(Poi(\Theta^{\smallsup{N}}_f)\neq Pois(\Theta_f))\\
    &=\prob( Poi(\Theta_f-\Theta^{\smallsup{N}}_f)\geq 1)\leq \expec[\Theta_f-\Theta^{\smallsup{N}}_f]
    =\expec\left[\sinv(U)-\sinv\Big(U_{\sss N}\Big)\right],\label{CL.Dist1}
    \end{align}
where the first inequality in the above chain is the coupling
inequality, the second the Markov inequality.

We first investigate the convergence of the integrals
$\expec\big[{(\Theta_f^{\smallsup{N}})^i}\big]$ to
$\expec\big[(\Theta_f)^i\big]$, for $i=1,2$, and then prove that this
is sufficient for conditions  (\mCo1)--(\mCo3).

\begin{lemma}[Convergence of moments of $\Theta_f^{\smallsup{N}}$]
\label{lem-conv-int} Let
$\Theta_f^{\smallsup{N}}=\sinv\Big(U_{\sss N}\Big)$, and
assume that $F$ satisfies \eqref{distribution}. Furthermore, assume
that $F$ is a continuous distribution function with density
$f\colon [0,\infty)\to [0,\infty)$, and that there exist an $a>0$
such that, as $y\downarrow 0$,
    \eq
    \label{f-cond}
    \int_y^1 \frac{\sinv(u)}{\underline{f}(\sinv(u))}\,\mathrm du =\mO{y^{-a}},
    \en
where $\underline{f}(x)=\inf_{0\leq y\leq x} f(y)$. Then, there
exists an $\alpha>0$ such that
    \eq
    \label{conv-int1}
    \int_0^1 \sinv(u)\br{\sinv(u)-\sinv(u_{\sss N})} \, \mathrm du =\mO{N^{-\alpha}},
    \en
and, for $i=1, 2$,
    \eq
    \label{conv-int2}
    \int_0^1 \sinv(u)^i-\sinv(u_{\sss N})^i \,\mathrm du =\mO{N^{-\alpha}}.
    \en
\end{lemma}

The proof of Lemma \ref{lem-conv-int} is followed by an example of cases where the conditions are satisfied.

\bigskip

{\noindent {\bf Proof of Lemma \ref{lem-conv-int}:}
We first claim that \eqref{distribution} implies that for $y\in (0,1)$,
        \begin{align}
        \sinv(y)&\leq c^b y^{-b},\label{CO2:claim1}
        \end{align}
where $b=1/(\tau-1)$.
We show the claim \eqref{CO2:claim1} by contradiction. Suppose that there is $y\in (0,1)$, satisfying
$ \sinv(y) > c^b y^{-b}$.
Observe that if  $\sinv(y)=w$, then ${\bar F}(x)>y$ for each $x<y$, by definition \eqref{invverd}.
Then with $x=c^b y^{-b}< \sinv(y)=w$,
\begin{equation}
{\bar F}(x)> y=cx^{1-\tau},
\end{equation}
which contradicts \eqref{distribution}. This proves \eqref{CO2:claim1}.

From the three statements in \eqref{conv-int1} and \eqref{conv-int2}, we only prove \eqref{conv-int2}
for $i=2$, the other two claims are similar, and in fact
easier to prove.
When $i=2$, we can rewrite and bound
    \begin{align}
    \int_0^1 \sinv(u)^2-\sinv(u_{\sss N})^2 \,\mathrm du
   \label{rew-bd-int}
    &\leq 2\int_0^1\sinv(u)\big(\sinv(u)-\sinv(u_{\sss N})\big)\,\mathrm du,
    \end{align}
where we have used the special product $x^2-y^2=(x-y)(x+y)$ and that  $\sinv(u_{\sss N})\leq \sinv(u)$,
for all $u\in (0,1)$. We split the integral into $u\in (0,N^{-d})$
and $u\in [N^{-d},1)$, for an appropriate $d>0$. For the former, we
bound using \eqref{CO2:claim1}
    \eq \nonumber
    \int_0^{N^{-d}}\sinv(u)\big(\sinv(u)-\sinv(u_{\sss N})\big)\,\mathrm du
    \leq \int_0^{N^{-d}}\sinv(u)^2 \,\mathrm du\leq c^{2b}\int_0^{N^{-d}}u^{-2b}\,\mathrm  du
    =\mO{N^{-d(1-2b)}},
    \en
where we used  that $b=1/(\tau-1)<1/2$, since $\tau>3$. For the integral over
$u\in [N^{-d},1)$, we note by Taylor's Theorem that
    \eq
    \sinv(u)-\sinv(u_{\sss N})
    =-(u-u_{\sss N})\frac{{\rm d}}{{\rm d} u}\sinv(u_{\sss N}^*)
    =\frac{u_{\sss N}-u}{f(\sinv(u_{\sss N}^*))},
    \en
for some $u_{\sss N}^*\in [u,u_{\sss N}]$. Since $u\mapsto
\sinv(u)$ is non-increasing, we have that $\sinv(u_{\sss N}^*)\leq
\sinv(u)$, so that
    \eq
    \frac{u_{\sss N}-u}{f(\sinv(u_{\sss N}^*))}\leq \frac{1}{N \underline{f}(\sinv(u))}.
    \en
Thus, we arrive at
    \eq
    \int_{N^{-d}}^1\sinv(u)\big(\sinv(u)-\sinv(u_{\sss N})\big)\,\mathrm du
    \leq \frac{1}{N}\int_{N^{-d}}^1 \frac{\sinv(u)}{\underline{f}(\sinv(u))}\,\mathrm du
    =\mO{N^{-1+ad}},
    \en
where, in the last inequality, we made use of \eqref{f-cond}. We
conclude that
    \eq
    \int_0^1 \sinv(u)^2-\sinv(u_{\sss N})^2 \,\mathrm du
    \leq \mO{N^{-d(1-2b)}}+\mO{N^{ad-1}} =\mO{N^{-\alpha}},
    \en
where $\alpha=\min\{d(1-2b), 1-ad\}>0$, provided we pick $d>0$ so
small that $1-ad>0$. \qed

\medskip

We now present the promised example:
\ch{Consider the class of those continuous distribution functions
$F$ satisfying \eqref{distribution} for which $f(x)=F'(x)$ is
positive and decreasing on $[0,\infty)$. We shall prove that such $F$
indeed satisfy \eqref{f-cond}.}
\ch{We will use \eqref{CO2:claim1} to prove \eqref{f-cond} when $f(x)=F'(x)$ is
positive and decreasing on $[0,\infty)$.}
\ch{Then, for such $f$,} $\underline{f}(x)=f(x)$. In this case, we
note that
    \eq
    2\frac{\sinv(u)}{f(\sinv(u))}
    =-\frac{{\rm d}}{{\rm d} u} \sinv(u)^2,
    \en
so that, by \eqref{CO2:claim1} with $b=1/(\tau-1)$,
    \eq
    2\int_y^1 \frac{\sinv(u)}{\underline{f}(\sinv(u))}\,\mathrm du
    =\sinv(y)^2\leq c^{2b} y^{-2b}=\mO{ y^{-2b}},
    \en
so that \eqref{f-cond} is satisfied.
In particular,
for
    \eq
    {\bar F}_{\sss \Lambda}(x)= x^{-\tau+1},\quad x\ge 1,
    \en
this proves the claim.
\vskip0.5cm

\bigskip

We will now show that the conditions (\mCo1)--(\mCo3) do hold for the
weights $\{\lambda_i\}_{i=1}^N=\sinv(i/N)$. We start with validating
the conditions (\mCo1) and (\mCo3), since for condition (\mCo2) more
effort is needed.

\bigskip

{\noindent \bf Condition (\mCo1):} Combining \eqref{muNnuN-def} with
\eqref{deflam} we rewrite $\mu_\sN$ as
    \eqn{\label{CL.mu.int}
    \mu_\sN=\frac{1}{N}\sum_{i=1}^N\bar
    F^{-1}(i/N) =\int_0^1\bar
    F^{-1}\br{u_\sN}\,\mathrm du,
    }
where as before $u_\sN=\frac{\bceil{Nu}}{N}$. Similarly, one can show that
    \eqn{\label{CL.nu.N}
    \nu_\sN=\frac{1}{\mu_\sN}\int_0^1
    \bar F^{-1}\br{u_\sN}^2\,\mathrm du. }
Therefore, we choose
    \begin{align}\label{CL.mu.nu}
    \mu&=\int_0^1\bar F^{-1}(u) \,\mathrm du& &\text{and}&
    \nu&=\frac{1}{\mu}\int_0^1\bar F^{-1}(u)^2  \,\mathrm du.&
    \end{align}
Effectively, \ch{we have replaced $u_\sN$ in \eqref{CL.mu.int}
and \eqref{CL.nu.N} by $u$.} Observe that the
following also hold
\begin{align}
    \mu&=\expec\big[\Theta_f\big],&
    \mu_\sN&=\expec\big[\Theta^{\smallsup{N}}_f\big],&
    \nu&=\expec\big[\Theta_f^2\big]/\expec\big[\Theta_f\big]&
    &\text{and}&
    \nu_\sN&=\expec\big[(\Theta^{\smallsup{N}}_f)^2\big]/\expec\big[\Theta^{\smallsup{N}}_f\big].&
\end{align}
Condition (\mCo1) is satisfied, since Lemma \ref{lem-conv-int},
where we take $i=1$,  implies the existence of a constant
$\alpha_1>0$ such that
    \eqn{\label{CL.afs.mu} |\mu_\sN-\mu|
    \leq \int_0^1\brc||{\bar F^{-1}\br{u} - \bar F^{-1}\br{u_\sN}}\,\mathrm du
    < N^{-\alpha_1}.
    }
Using \eqref{CL.nu.N} and \eqref{CL.mu.nu}, we obtain that
    \begin{align*}
    |\nu_\sN-\nu| &= \brc||{\frac{1}{\mu_\sN}\int_0^1\bar
    F^{-1}(u_\sN)^2  \,\mathrm du-\frac{1}{\mu}\int_0^1\bar
    F^{-1}(u)^2  \,\mathrm du}
    \\&\leq
    \frac{1}{\mu_\sN}\int_0^1\brc||{\bar
    F^{-1}\br{u_\sN}^2-F^{-1}\br{u_\sN}^2}\,\mathrm du
    +\br{\frac{1}{\mu_\sN}-\frac{1}{\mu}}\int_0^1\bar
    F^{-1}\br{u}^2\,\mathrm du.
    \end{align*}
We bound the above display, using Lemma \ref{lem-conv-int} and
\eqref{CL.afs.mu}, by
    $$
    |\nu_\sN-\nu| \leq \frac{1}{\mu-N^{-\alpha_1}}N^{-\alpha_1}
    +\br{\frac{1}{\mu-N^{-\alpha_1}}-\frac{1}{\mu}}\int_0^1\bar
    F^{-1}\br{u}^2\,\mathrm du =\mO{N^{-\alpha_1}},
    $$
since
    \eqn{\label{CL.2m} \int_0^1\bar F^{-1}\br{u}^2\,\mathrm du
    =\mexpec{\bar  F^{-1}\br{U}^2}=\mexpec{\Theta_f^2} < \infty,
    }
where $\Theta_f=\bar F^{-1}\br{U}$ is a random variable with
distribution $F$ satisfying \eqref{distribution}.

\bigskip

{\noindent \bf Condition (\mCo3):}
Observe that, \ch{again using that $u\mapsto \sinv(u)$ is non-increasing,}
    \begin{align}
    \frac{1}{N}\sum_{i=1}^N \ml i^{\tau-1-\varepsilon}\nonumber
    &=\frac{1}{N}\sum_{i=1}^N \bar F^{-1}(i/N)^{\tau-1-\varepsilon}
    =\mexpec{\bar F^{-1}(U_{\sss N})^{\tau-1-\varepsilon}}
    \leq \mexpec{\bar F^{-1}(U)^{\tau-1-\varepsilon}},
    \end{align}
and $\bar F^{-1}(U)$  is a random variable with distribution $F$
which satisfies \eqref{distribution}, compare with \eqref{CL.2m}.
Therefore $\mexpec{\bar F^{-1}(U)^{\tau-1-\varepsilon}}<\infty$,
which in turn implies that
$$
\limsup_{N\rightarrow \infty} \frac{1}{N}\sum_{i=1}^N \ml
i^{\tau-1-\varepsilon}
 < \infty.
$$

Since, for fixed $N$, we have that $\bar F^{-1}(i/N)$ is non-increasing in
$i$,  we obtain from \eqref{CO2:claim1},
$$
\lambda_\sN^\smallsup N =\max_{1
\leq i \leq N} \bar F^{-1}(i/N)=\bar F^{-1}(1/N)=\mO{N^{1/(1-\tau)}} \leq
N^{\gamma}.
$$

\bigskip

{\noindent \bf Condition (\mCo2):} We bound \eqref{CL.Dist1} by
Lemma \ref{lem-conv-int}: there exists a constant $\alpha_2>0$
such that
$$
\mbox{d}_{\sss \rm
TV}(Poi(\Theta^{\smallsup{N}}_f),Poi(\Theta_f))\leq
\expec\left[\sinv(U)-\sinv\Big(U_{\sss N}\Big)\right]=\mO{N^{-\alpha_2}}.
$$
Similarly as \eqref{CL.Dist1}, we obtain
    \begin{align}\nonumber
    \frac12 \sum_{n=0}^\infty |\mg n-\mgl n|&= \mbox{d}_{\sss \rm TV}(\mPoisson{\Theta^{\smallsup{N}}_g},\mPoisson{\Theta_g})\\
    &\leq \expec[|\Theta_g-\Theta^{\smallsup{N}}_g|]
    =\mexpec{\brc||{ \bar F^{-1}_{\Theta_g}(U)-\bar F^{-1}_{\Theta_g^{\smallsup{N}}}(U)}}.
\label{CL.Dist2}
    \end{align}

\ch{We claim that for two distribution functions
$G\colon \mathbb R \mapsto [0,1]$ and $H\colon \mathbb R \mapsto [0,1]$
with $G(0)=H(0)=0$, we have
    \eqn{\label{L1-stat}
     \ch{\int_0^1}\brc||{\bar G^{-1}(u)-\bar H^{-1}(u)}\,\mathrm{d}u
    =\int_0^\infty\brc||{\bar G(x)-\bar H(x)}\,\mathrm{d}x.
    }
To see this, note that when $\bar G(x)\geq \bar H(x)$ for every $x\in {\mathbb R}$,
so that $\bar G(x)-\bar H(x)\geq 0, \bar G^{-1}(u)-\bar H^{-1}(u)\geq 0$,
then both sides of \eqref{L1-stat} are equal to $\expec[X]-\expec[Y]$,
where $X$ and $Y$ have distribution functions $G$ and $H$ respectively.
This proves the claim for $G$ and $H$ for which $\bar G(x)\geq \bar H(x)$ for
every $x\in {\mathbb R}$. For other $G$ and $H$, we write,
for $x\in [0,\infty)$,
    \eqn{
    (\bar G\vee \bar H)(x)=\max\{\bar G(x), \bar H(x)\},
    \qquad
    (\bar G\wedge \bar H)(x)=\min\{\bar G(x), \bar H(x)\},
    }
so that
    \eqn{
    \int_0^\infty \brc||{\bar G(x)-\bar H(x)}\,\mathrm{d}x
    =\int_0^\infty(\bar G\vee \bar H)(x)-(\bar G\wedge \bar H)(x)\,\mathrm{d}x.
    }
Then, both $\bar G\wedge \bar H$ and $\bar G\vee \bar H$ are
survival functions for which $(\bar G\vee \bar H)(x)\geq (\bar G\wedge \bar H)(x)$
for every $x\in {\mathbb R}$. Thus, by the above proof, \eqref{L1-stat}
holds for them. The claim is then completed by noting that
    \eqn{
    (\bar G\vee \bar H)^{-1}(u)=\inf \{s: \max\{\bar G(s), \bar H(s)\}\leq u\}
    =\max\{\bar G^{-1}(u), \bar H^{-1}(u)\},
    }
and
    \eqn{
    (\bar G\wedge \bar H)^{-1}(u)=\inf\{s: \min\{\bar G(s), \bar H(s)\}\leq u\}
    =\min\{\bar G^{-1}(u), \bar H^{-1}(u)\}
    }
so that
    \eqn{
    (\bar G\vee \bar H)^{-1}(u)-(\bar G\wedge \bar H)^{-1}(u)
    =\big|\bar G^{-1}(u)-\bar H^{-1}(u)\big|.
    }
This completes the proof of \eqref{L1-stat}.
}

Using \eqref{L1-stat}, we can rewrite display \eqref{CL.Dist2} as
    \eqn{\label{CL.dit} \mbox{d}_{\sss \rm
    TV}(Poi(\Theta^{\smallsup{N}}_g),Poi(\Theta_g)) \leq
    \int_{0}^1{\brc||{ \bar F^{-1}_{\Theta_g}(u)-\bar F^{-1}_{\Theta_g^{\smallsup{N}}}(u)}}\,\mathrm d u=
        \int_{0}^\infty \brc||{\bar F_{\Theta_g}(x)-\bar F_{\Theta_g^{\smallsup{N}}}(x)} \,\mathrm d x.
    }
Recall that $\mu=\mexpec{\Theta_f}$ and
$\mu_\sN=\mexpec{\Theta^{\smallsup{N}}_f}$, then, using Lemma
\ref{lem.mixingpoisson}, the integrand can be rewritten as
    \begin{align}
    \brc||{\bar F_{\Theta_g}(x)-\bar F_{\Theta_g^{\smallsup{N}}}(x)}
    &
    =\brc||{\mexpec{{\frac{\Theta_f}{\mu}\indic{\Theta_f>x}-\frac{\Theta_f^{\smallsup{N}}}{\mu_\sN}\indic{\Theta_f^{\smallsup{N}}>x}}}}
    \nonumber\\&\leq \mexpec{\brc||{\frac{\bar F^{-1}(U)}{\mu}\indic{\bar F^{-1}(U)>x}-\frac{\bar  F^{-1}(U_{\sN})}{\mu_\sN}
    \indic{\bar F^{-1}(U_{\sN}>x}}}
    \nonumber\\&=\int_{0}^1 \brc||{\frac{\bar F^{-1}(u)}{\mu}\indic{\bar F^{-1}(u)>x}-\frac{\bar  F^{-1}(u_\sN)}{\mu_\sN}\indic{\bar F^{-1}(u_\sN)>x}}\,\mathrm du.\label{CL.dit.int}
    \end{align}
Using \eqref{CL.dit},  \eqref{CL.dit.int}, and the triangle
inequality twice, we obtain
\begin{align}
\mbox{d}_{\sss \rm TV}(Poi(\Theta^{\smallsup{N}}_g),Poi(\Theta_g))
    &\leq \brc||{\frac{1}{\mu}-\frac1{\mu_\sN}}\int_{0}^\infty\int_{0}^1 \bar F^{-1}(u)\indic{\bar F^{-1}(u)>x}\,\mathrm du \,\mathrm dx\nonumber
    \\&\quad+\frac{1}{\mu_\sN}\int_{0}^\infty\int_{0}^1
    \bar F^{-1}(u)\brc||{\indic{\bar F^{-1}(u)>x}-\indic{\bar F^{-1}(u_\sN)>x}}\,\mathrm du \,\mathrm dx \nonumber
    \\&\quad+\frac{1}{\mu_\sN}\int_{0}^\infty\int_{0}^1
    \brc||{\bar F^{-1}(u)-\bar F^{-1}(u_\sN)}\indic{\bar F^{-1}(u_\sN)>x}\,\mathrm du \,\mathrm dx. \label{CL.di}
\end{align}

Next, we will bound each double integral in display \eqref{CL.di}. The
first double integral is finite, since
    \eqn{\label{CL.dd.1}
    \int_{0}^\infty\int_{0}^1 \bar F^{-1}(u)\indic{\bar
    F^{-1}(u)>x}\,\mathrm du \,\mathrm dx
    =\int_{0}^1 \bar F^{-1}(u)^2 \,\mathrm du=\mexpec{\Theta_f^2}<\infty,
    }
where we used in the last step that the distribution $F$ satisfies
\eqref{distribution}. Observe that $\indic{\bar F^{-1}(u)>x}\not
=\indic{\bar F^{-1}(u_\sN)>x}$ if and only if $\bar F^{-1}(u_\sN)<
x < \bar F^{-1}(u)$, thus
    \eqn{\label{CL.dd.2}
    \int_{0}^\infty\int_{u=0}^1
    \bar F^{-1}(u)\brc||{\indic{\bar F^{-1}(u)>x}-\indic{\bar F^{-1}(u_\sN)>x}}\,\mathrm du \,\mathrm dx
    =\int_{0}^1 \bar F^{-1}(u)\br{\bar F^{-1}(u)- \bar F^{-1}(u_\sN)}\,\mathrm du.
    }
Using \eqref{conv-int1} it follows that the right hand side of \eqref{CL.dd.2} is $\mO{N^{-\alpha}}$.
Finally, using again \eqref{conv-int1},
\begin{multline}\label{CL.dd.3}
\int_{0}^\infty\int_{0}^1
    \br{\bar F^{-1}(u)-\bar F^{-1}(u_\sN)}\indic{\bar F^{-1}(u_\sN)>x}\,\mathrm du \,\mathrm dx
    =\int_{0}^1
    \br{\bar F^{-1}(u)-\bar F^{-1}(u_\sN)}\bar F^{-1}(u_\sN)\,\mathrm du
    \\ \leq  \int_{0}^1
    \bar F^{-1}(u)\br{\bar F^{-1}(u)-\bar F^{-1}(u_\sN)}\,\mathrm du
    =\mO{N^{-\alpha}}.
\end{multline}

    Combining \eqref{CL.dd.1}--\eqref{CL.dd.3} with \eqref{CL.di}, we obtain that
        $$
        \mbox{d}_{\sss \rm TV}(Poi(\Theta^{\smallsup{N}}_g),Poi(\Theta_g))
        =\mO{\frac{1}{\mu}-\frac1{\mu_\sN}}+\frac{1}{\mu_\sN}\mO{N^{-\alpha}}+\frac{1}{\mu_\sN}\mO{N^{-\alpha}}=\mO{N^{-\alpha_3}},
        $$
    where we used that $|\mu-\mu_\sN|=\mO{N^{-\alpha_1}}$ and take $\alpha_3=\min\{\alpha,\alpha_1\}$.

\newcommand{\msexpec}[2][\sN]{\mathbb E_{#1}\hspace{-.3em}\brh{#2}}
\newcommand{\msprob}[2][\sN]{\mathbb P_{\hspace{-.15em} #1} \hspace{-.20em} \br{#2}}

\section{Appendix: i.i.d.\ capacities} \label{iidcase}
\ch{In this section, we prove Theorem \ref{corIIDcase} by
showing that the results of the Poissonian random
graph with deterministic capacities also hold when the capacities are
replaced by i.i.d.\ capacities satisfying \eqref{distribution}.}

\ch{More precisely, we associate a sequence $\{\mL i\}_{i=1}^N$
of positive i.i.d.\ random variables to the nodes,
where the random capacities $\{\mL i\}_{i=1}^N$  have distribution
function $F_\Lambda(x)=\mprob{\mL\nill \leq x}$, and we set
the capacity of node $i$ equal to $\mL i$, for $1\leq i \leq N$.}

Define
  \eqn{\label{defJ}
  \ch{{\cal J}=}{\cal J}(\alpha)={\cal J}_1(\alpha) \cap {\cal J}_2(\alpha) \cap {\cal J}_3,
  }
where
\begin{align}
    {\cal J}_1(\alpha)&=
    \bra{|\mu_\sN-\mu|  < N^{-\alpha}}
    \cap
    \bra{|\nu_\sN -\nu| < N^{-\alpha}},\\
    {\cal J}_2(\alpha)&=
            \bra{ d_{\sss \rm TV}(\mf\nill, \mfl\nill) < N^{-\alpha}}
    \cap
            \bra{ d_{\sss \rm TV}(\mg\nill, \mgl\nill) < N^{-\alpha}},
    \intertext{and, \ch{for some $B>0$ sufficiently large,}}
        \label{defJ3}
    {\cal J}_3&=\bra{\sum_{n=0}^{\infty}
    n^{\tau-2-\vep}\mg n\ch{\leq B}}
    \cap \bra{\Lambda_\sN^\smallsup N \leq N^{\gamma}},
 \end{align}
and where the \ch{exponent} $\gamma$ is equal to $1/(\tau-1)+\varepsilon$
\ch{(recall \eqref{P:assumP}).
 }

In Sections \ref{convofmean}--\ref{sec-iidCo3} below, we will
show that there exist constants $$\alpha_1, \alpha_2, \beta_1, \beta_2, \beta_3>0$$
such that  the events ${\cal J}_1(\alpha_1)$, ${\cal J}_2(\alpha_2)$, ${\cal J}_3$
occur with probability exceeding $1-\mO{N^{-\beta_1}}$, $1-\mO{N^{-\beta_2}}$ and
$1-\mO{N^{-\beta_3}}$, respectively.
Then, by taking $\alpha=\min\{\alpha_1,\alpha_2\}$ and
$\beta=\min\{\beta_1,\beta_2,\beta_3\}$ the event ${\cal J}(\alpha)$ occurs
with probability exceeding $1-\mO{N^{-\beta}}$.
On the event ${\cal J}(\alpha)$ the conditions (\mCo1)--(\mCo3) do hold,
and this then proves Theorem \ref{corIIDcase}.

For the sake of completeness, we give the counterparts of
    \eqref{muNnuN-def} and \eqref{def-fgn-n} using the random capacities.
    Display \eqref{muNnuN-def} becomes
    \begin{align}
   \label{muNnuN-def-rand}
    \mu_{\sN}&=\frac 1N \sum_{i=1}^N \mL i&
    &\text{and}&
    \nu_{\sN}&=\frac{1}{\mu_\sN N}\sum_{i=1}^N \mL i^2,&
    \end{align}
Using the law of large number the values $\mu_\sN$ and $\nu_\sN$ converge a.s.\ to
\ch{$\mu$ and $\nu$ in \eqref{pick idd}, respectively.}

The offspring distribution of the first and second generation is given by
(compare \eqref{def-fgn-n}),
\begin{align}
    \label{def-fgn-n-rand}
    \mf n
    &=\frac{1}{N}\sum_{i=1}^N e^{-\mL i}\frac{\mL i^n}{n!}&
    &\text{and}&
    \mg n&=\frac{1}{\mu_{\sN}N}\sum_{i=1}^N e^{-\mL i}\frac{\mL i^{n+1}}{n!}
    =\frac{(n+1)\mf {n+1}}{\mu_\sN},&
    \end{align}
which, by the strong law of large numbers, converge a.s.\ to $f$ and $g$ in
\eqref{pick idd}, respectively.

\subsection{Convergence of means (\mCo1):}
\label{convofmean}

In this section we will show that condition (\mCo1) holds
\whp, i.e, that there exist constants $\alpha_1,\beta_1>0$ such
that the event
    $$
    {\cal J}_1(\alpha_1)=\bra{|\mu_\sN-\mu| < N^{-\alpha_1}} \cap \bra{|\nu_\sN -\nu| < N^{-\alpha_1}}
    $$
    occurs with probability exceeding $1-\mO{N^{-\beta_1}}$,
    where $\mu$ and $\nu$ are \ch{defined in \eqref{pick idd}.}

The next lemma is the crucial \ch{estimate} of the proof:
\begin{lemma} \label{lem:SnA} Let $q \in (0, \tau-1)$ and define
    $$
    S_{\sN,q}=\frac{1}{N}\sum_{i=1}^N \mL i^q
    \qquad \text{and}
    \qquad
    {\cal S}_{\sN,q}(u)=\bra{\brc||{S_{\sN,q} - \mexpec{\mL\nill^q}} \leq  N^{-u}},
    $$
then there exist constants $u, v>0$
such that
\begin{align}\label{claim S:A}
\mprob{{\cal S}_{\sN,q}(u)} =1-\mO{N^{-v}}.
\end{align}
\end{lemma}

{\noindent \bf Proof.} Apply the Marcienkiewicz-Zygmund inequality \cite{Gut1},
with $X_i=\Lambda_i^q-\expec[\Lambda_i^q]$ and $r\in (1,2)$, such that $qr<\tau-1$.
Then,
\begin{eqnarray*}
\mprob{{\cal S}_{\sN,q}(u)^c}
&=&\mprob{\Big|\frac1{N} \sum_{i=1}^N X_i\Big| >  N^{-u}}
=\mprob{\Big|\frac1{N} \sum_{i=1}^N X_i\Big|^r >  N^{-ur}}\\
&= & \mprob{\Big| \sum_{i=1}^N X_i\Big|^r >  N^{r(1-u)}}\le N^{-r(1-u)}c_r N\expec[|X_1|^r]
={\cal O}(N^{1-r(1-u)}),
\end{eqnarray*}
\ch{since} $\expec[|X_1|^r]\leq \expec[|\Lambda_1|^{qr}]<\infty $.
\ch{Now recall that $r>1$ and fix $u>0$ so small that $v=r(1-u)-1>0$.}
\qed

\medskip

    Using Lemma \ref{lem:SnA}, we pick the constants $u,v>0$ such that
    $$
    \mprob{{\cal S}_{\sN,1}(u) \cap {\cal S}_{\sN,2}(u)} =1-\mO{N^{-v}}.
    $$
    Observe that $\mu_{\sN}=S_{\sN,1}$ and $\nu_{\sN}=S_{\sN,2}/S_{\sN,1}$.
    On the event ${\cal S}_{\sN,1}(u) \cap {\cal S}_{\sN,2}(u)$, we have that
    $$
    |\mu_\sN-\mu|=|S_{\sN,1}-\mexpec{\mL\nill}| =\mO{ N^{-u}},
    $$
    and
    \begin{align}\nonumber
    |\nu_\sN-\nu|&=\brc||{\frac{S_{\sN,2}}{S_{\sN,1}}-\frac{\mexpec{\mL\nill^2}}{\mexpec{\mL\nill}}}
    \leq S_{\sN,2}\brc||{\frac{1}{S_{\sN,1}}-\frac{1}{\mexpec{\mL\nill}}} + \frac{|S_{\sN,2}-\mexpec{\ch{\Lambda^2}}|}{\mexpec{\mL\nill}}
    =\frac{S_{\sN,2}|\mexpec{\mL\nill}-S_{\sN,1}|}{\mexpec{\mL\nill} S_{\sN,1}} +\mO{N^{-u}}    \\&=\frac{\br{\mexpec{\mL\nill^2}+\mO{N^{-u}}}\mO{N^{-u}}}{\mexpec{\mL\nill}(\mexpec{\mL\nill}-\mO{N^{-u}})} +\mO{N^{-u}}
    =\mO{N^{-u}}.
    \end{align}
    Therefore, we can pick $\alpha_1=u$ and $\beta_1=v$.
    \qed

    \bigskip

\subsection{Convergence of branching processes (\mCo2):}
\label{sec-iidCo2}
    In this section we will show that there exist constants $\alpha_2,\beta_2>0$ such that
    \begin{equation}
    \label{totvarg}
    \mprob{{\rm d}_{\sss \rm TV}(\mg\nill, \mgl\nill) > N^{-\alpha_2} } =\mO{N^{-\beta_2}}.
    \end{equation}
    The proof of the existence of constants $\alpha_2,\beta_2>0$ such that
    $\mprob{{\rm d}_{\sss \rm TV}(\mf\nill, \mfl\nill) > N^{-\alpha_2} } =\mO{N^{-\beta_2}}$ is
    less difficult and therefore omitted.

    In this section, we add a subscript $N$ to $\mathbb P$ and $\mathbb E$,
    to make clear that we condition on the \ch{i.i.d.} capacities $\{\mL i\}_{i=1}^N$. Thus,
    \begin{align}
    \msprob{\,\cdot\,}&=\mprob{\, \cdot \,|\mL 1, \ldots, \mL N}&
        &\text{and}&
        \msexpec{\,\cdot\,}&=\mexpec{\, \cdot \,|\mL 1, \ldots, \mL N}.&
    \end{align}

    Using Lemma \ref{lem.mixingpoisson} with $\{\mf n\}_{n\geq 0}$ and $\{\mfl n\}_{n\geq 0}$, it follows that
    $
    \mg n=\msprob{\mPoisson{\Theta^\smallsup{N}_g}=n},
    $
    where, for $i=1,2,\ldots,N$, $$\msprob{\Theta^{\smallsup{N}}_g=\mL i}=\frac{\mL i}{N\mu_\sN},$$
        and $\mgl n=\mprob{\mPoisson{\Theta_g}=n}$,
    where $\Theta_g$ is given by
    \eqn{\label{def:FW}
    F_{\Theta_g}(x)=\frac{1}{\mu}\mexpec{\mL\nill\indic{\mL\nill \leq x}}.
    }
    Let $U$ uniformly distributed over $(0,1)$, then we couple $ \Theta^\smallsup{N}_g$ and $\mPoisson{\Theta_g}$ by letting
    $\Theta_g=\ch{\sinv}_{\Theta_g}(U)$ and $\Theta^\smallsup{N}_g=\ch{\sinv}_{\Theta^\smallsup{N}_g}(U)$.
    Arguing as in \eqref{CL.Dist1}, we find that,
    \begin{align*}
    {\rm d}_{\sss \rm TV}(\mg\nill, \mgl\nill)\nonumber
    &={\rm d}_{\sss \rm TV}(\mPoisson{\Theta^\smallsup{N}_g},\mPoisson{\Theta_g})
    \leq
     \msexpec{\brc||{\Theta^\smallsup{N}_g-\Theta_g}}.
    \end{align*}
Then, using the coupling \ch{and \eqref{L1-stat}, we can write}
    \eqan{\label{totalvar}
    {\rm d}_{\sss \rm TV}(\mg\nill, \mgl\nill)
        &\leq\msexpec{\brc||{\Theta^\smallsup{N}_g-\Theta_g}}
        =\msexpec{\brc||{\ch{\sinv}_{\Theta^{\smallsup{N}}_g}(U)-\ch{\sinv}_{\Theta_g}(U)}}\nn\\
        &=\int_{0}^1\brc||{\ch{\sinv}_{\Theta^{\smallsup{N}}_g}(u)-\ch{\sinv}_{\Theta_g}(u)} \,\mathrm du
        \ch{=\int_{0}^{\infty}\brc||{\bar{F}_{\Theta^{\smallsup{N}}_g}(x)-\bar{F}_{\Theta_g}(x)} \,\mathrm dx,}
    }
    where $F_{{\Theta_g}}$ is given by \eqref{def:FW} and $F_{{\Theta^{\smallsup{N}}_g}}$ by
    \eqn{\label{def:FWN}
    F_{{\Theta^{\smallsup{N}}_g}}(x)=\frac1{\mu_\sN N}\sum_{i=1}^N \mL i\indic{\mL i\leq x}.
    }

    Before we prove \eqref{totvarg}, we state two claims. The proofs of these claims are deferred to the end of this section.
    The first claim is that
    \eqn{ \label{B2:cond1}
    \ch{{\rm d}_{\sss \rm TV}(\mg\nill, \mgl\nill)\leq}\brc||{1-\frac{\mu_\sN}\mu}\nu_\sN  +
    \frac{E_\sN}{\mu},
    }
    where
    \eqn
    {\label{G:E}
    E_\sN=\int_{0}^\infty \brc||{\frac{1}{N}\sum_{i=1}^N \mL i \indic{\mL i >x} -\mexpec{\mL\nill\indic{\mL\nill > x}}} \,\mathrm d x .
    }
    Secondly, we claim that there exist constants $\alpha,\beta>0$ such that
    \eqn{\label{claimE}
    \mprob{E_\sN > N^{-\alpha}} =\mO{N^{-\beta}}.
    }

    We finish the proof of \eqref{totvarg}, given the claims \eqref{B2:cond1} and \eqref{claimE}.
    Indeed, from Section \ref{convofmean} we know that there exist constants $\alpha_1,\beta_1>0$ such that $\mprob{{\cal J}_1(\alpha_1)} \geq 1-\mO{N^{-\beta_1}}$.

    Combining \eqref{B2:cond1} and \eqref{claimE}, we can bound the total variation by
    \begin{align}
    \mprob{{\rm d}_{\sss \rm TV}(\mg\nill, \mgl\nill) > N^{-\alpha_2} }
    \nonumber&\leq N^{\alpha_2}\mexpec{\br{\brc||{1-\frac{\mu_\sN}\mu}\nu_\sN  + E_\sN/\mu} \mindic{ \bra{E_\sN \leq N^{-\alpha}}\cap {\cal J}_1(\alpha_1)}}
    \nonumber\\&\quad\quad+\mprob{E_\sN>N^{-\alpha}} + \mO{N^{-\beta_1}}
    =\mO{N^{-\min\{\alpha_1-\alpha_2,\alpha-\alpha_2,\beta,\beta_1\}}}.
    \end{align}
    Pick $\alpha_2 \in(0, \min\{\alpha_1,\alpha\})$, then the right-hand side in the above display tends to zero as $N$ tends to
    infinity, and \eqref{totvarg} follows, when \eqref{B2:cond1} and
    \eqref{claimE} hold.
\qed

\bigskip

Finally, we finish the proof by showing, in turn, the claims \eqref{B2:cond1} and \eqref{claimE}.
\bigskip

    {\noindent \bf Proof of \eqref{B2:cond1}:}
    Using \eqref{def:FW}, \eqref{totalvar} and \eqref{def:FWN} we bound $\ch{{\rm d}_{\sss \rm TV}(\mg\nill, \mgl\nill)}$ by
    \begin{align}
    {\rm d}_{\sss \rm TV}(\mg\nill, \mgl\nill)
    \nonumber&
    = \int_0^\infty \brc||{\frac1{N\mu_\sN}\sum_{i=1}^N \mL i \indic{\mL i > x}-\frac{1}{\mu}\mexpec{\mL\nill\indic{\mL\nill > x}}} \mathrm d x
    \nonumber\\&\leq
    \brc||{1-\frac{\mu_\sN}\mu}\int_0^\infty{\frac1{N\mu_\sN}\sum_{i=1}^N \mL i \indic{\mL i > x}} \mathrm dx
    \nonumber\\&\quad\quad+
    \frac1\mu\int_0^\infty\brc||{ \frac1{N}\sum_{i=1}^N \mL i \indic{\mL i > x}-\mexpec{\mL\nill\indic{\mL\nill >x}}} \mathrm d x
    \nonumber\\&=\brc||{1-\frac{\mu_\sN}\mu}\nu_\sN+ E_\sN/\mu,
    \end{align}
    where in the last step we used that
    $$
    \int_0^\infty{\frac1{N\mu_\sN}\sum_{i=1}^N \mL i \indic{\mL i > x}} \mathrm d x =\frac1{N\mu_\sN}\sum_{i=1}^N \mL i
    \int_0^\infty{ \indic{\mL i > x}} \mathrm d x
    =\frac{1}{N\mu_\sN}\sum_{i=1}^N \mL i^2=\nu_\sN.
    $$
    \ch{This proves \eqref{B2:cond1}.}\qed

    \bigskip

{\noindent \bf Proof of \eqref{claimE}:} \ch{For $a>1$,
we define the density $f_a$ by
    \eqn{
    f_a(x)=\frac{a-1}{(x+1)^a}.
    }
Then, for any $\delta, a>1$, by}
subsequently using  the Markov inequality and then Jensen's inequality,
\ch{
    \begin{align}
    \mprob{E_\sN > N^{-\alpha}}
    & =\mprob{E_\sN^{\delta} > N^{-\alpha\delta}}\label{EN:als:int}\\
    &\leq
    N^{\alpha\delta}(a-1)^{-\delta}\mexpec{\br{\int_{0}^\infty (1+x)^a\brc||{\frac{1}{N}\sum_{i=1}^N \mL i \indic{\mL i >x} -
    \mexpec{\mL\nill\indic{\mL\nill > x}}}f_a(x)    \,\mathrm d x }^\delta}\nonumber
    \\&\leq  N^{\alpha\delta}\mexpec{{\int_{0}^\infty (1+x)^{a(\delta-1)}\brc||{\frac{1}{N}\sum_{i=1}^N \mL i \indic{\mL i >x}
     -\mexpec{\mL\nill\indic{\mL\nill > x}}}^\delta
    \,\mathrm d x }}\nonumber
    \\&=  N^{\alpha\delta}\int_{0}^\infty(1+x)^{a(\delta-1)}
    \mexpec{{ \brc||{\frac{1}{N}\sum_{i=1}^N \mL i \indic{\mL i >x} -\mexpec{\mL\nill\indic{\mL\nill > x}}}^\delta
   }} \,\mathrm d x
    .\nonumber
    \end{align}
In order to bound \eqref{EN:als:int}, we use the
Marcinkiewicz-Zygmund inequality \cite[Chapter 3; Corollary
8.2]{Gut1}
    on the variables  $X_i(x)=\Lambda_i \indic{\mL i >x} $, $1 \leq i \leq N$.
    Then we can bound \eqref{EN:als:int} from above by
    \begin{align}\nonumber
    \mprob{E_\sN > N^{-\alpha}}
    &\leq {\cal O}(N^{\alpha\delta}) \int_0^\infty (1+x)^{a(\delta-1)}\mexpec{\brc||{\frac{1}{N}\sum_{i=1}^N X_i(x) -\mexpec{X(x)}}^\delta}
    \,\mathrm d x \\
    &\leq {\cal O}(N^{\alpha\delta+1-\delta})\int_0^\infty (1+x)^{a(\delta-1)}
    \mexpec{\brc||{X(x)-\mexpec{X(x)}}^\delta} \,\mathrm d x.
    \label{fin0}
    \end{align}
    Since for any $x>0$,
    \begin{align}
    \mexpec{\brc||{X(x)-\mexpec{X(x)}}^\delta}
        \leq \mexpec {X(x)^{\delta}}
        =\int_x^\infty y^{\delta}\, \mathrm d F_{\Lambda}(y),
                    \end{align}
we obtain
\eqan
{ \int_0^\infty (1+x)^{a(\delta-1)}\mexpec{\brc||{X(x)-\mexpec{X(x)}}^\delta} \,\mathrm d x
&\le\int_{x=0}^{\infty} (1+x)^{a(\delta-1)}\int_{y=x}^\infty y^{\delta}\, \mathrm d F_{\Lambda}(y)\,\mathrm d x\nonumber\\
&\leq \int_0^\infty y^{\delta}(1+y)^{a(\delta-1)+1} \, \mathrm d F_{\Lambda}(y)<\infty,
\label{fin1}
}
when we first pick $a>1$ and next $\delta>1$ so small that $\delta+a(\delta-1)<\tau-2$.

Combining \eqref{fin0} and \eqref{fin1}, we obtain
    \eqn{\label{doel:deltaalpha}
    \mprob{E_\sN> N^{-\alpha}}
       =\mO{N^{\alpha\delta+1-\delta}}=\mO{N^{-(\delta(1-\alpha)-1)}},
    }
The right hand side of \eqref{doel:deltaalpha}, tends to zero if $\delta(1-\alpha)-1>0$,
which is the case if we choose $\alpha$ small and $\delta$ slightly bigger than $1$.
}\qed

\subsection{Moment and maximal bound on capacities (\mCo3):}
\label{sec-iidCo3}
In this subsection we show that
$\mprob{{\cal J}_3}=1-\mO{N^{-\beta_3}}$.
Observe, that for every  $\xi >0$, using Boole's inequality,
    $$
    \mprob{\mL \sN ^\smallsup{\sN} \geq N^{\xi}} \leq
    \sum_{i=1}^N \mprob{\mL i \geq N^{\xi}}
    =\mO{N^{1-(\tau-1)\xi}}=\mO{N^{-((\tau-1)\xi-1)}}.
    $$
If $\xi>1/(\tau-1)$ then the right hand side in the above display tends to zero.  Remember that $\gamma=1/(\tau-1)+\varepsilon$, therefore we have
    $$
    \mprob{\mL \sN ^\smallsup{\sN} \geq N^{\gamma}} =\mO{N^{-(\tau-1)\varepsilon}}.
    $$
\ch{This proves the bound on the maximal capacity in (\mCo3).}

Apply Lemma \ref{lem:SnA} with $q=\tau-1-\varepsilon$, which states that there exist constants $u,v>0$ such that
    $$
    1-\mprob{{\cal S}_{\sN,\tau-1-\varepsilon}(u)}=\mprob{\brc||{\frac 1N \sum_{i=1}^{N} \mL i^{\tau-1-\varepsilon} - \mexpec{\mL\nill^{\tau-1-\varepsilon}}} > N^{-u} } = \mO{N^{-v}}.
    $$
Furthermore, \ch{\eqref{helper:mom} in the proof in
Lemma \ref{lem-gvslambda} shows that
there exist constants $C_1, C_2<\infty$ such that,
uniformly in $\{\Lambda_i\}_{i=1}^N$,
    \eqn{
    \sum_{n=0}^{\infty}
    n^{\tau-2-\vep}\mg n\leq C_1+C_2\frac 1N \sum_{i=1}^{N} \mL i^{\tau-1-\varepsilon}.
    }
    }
Therefore, \ch{for $N$ sufficiently large and  $B>C_1+C_2\mexpec{\mL\nill^{\tau-1-\varepsilon}}$,}
    \begin{align}
    1-\mprob{{\cal J}_3} \nonumber
    &\leq \mprob{\mL \sN ^\smallsup{\sN} \geq N^{\gamma}}+
    \mprob{\bra{\ch{\sum_{n=0}^{\infty}
    n^{\tau-2-\vep}\mg n>B}} \cap {\cal S}_{\sN,\tau-1-\varepsilon}(u) }+\mO{N^{-v}}\nonumber
    \\&\leq \mO{N^{-(\tau-1)\varepsilon}}+0+\mO{N^{-v}}
    =\mO{N^{-\min\{(\tau-1)\varepsilon,v\}}}.
    \end{align}
Thus, \ch{we may} take $\beta_3=\min\{(\tau-1)\varepsilon,u\}$.

\switchText{}{\section{Appendix: Coupling of shortest path graphs to branching processes}\label{appendix coupling}
\ch{In this appendix, we prove Proposition \ref{coupling of NR} in
Section \ref{coupling of NR:A} and Proposition \ref{coupling of sums}
in Section \ref{proposition:coupling of sums:bewijs}.}

\subsection{Proof of Proposition \ref{coupling of NR}}
\label{coupling of NR:A}
In this part of the appendix we give the main result on the coupling between the
capacities of the \mRNprocess\ and those of the \mrRNprocess.
For convenience we restate Proposition \ref{coupling of NR} as:
\begin{prop} \label{coupling of NR:B}
    There exist constants $\eta, u_1,v_1>0$, such that for
    all $t \leq (1/2+\eta)\log_\nu N$,
    \begin{align}
    \label{mainprop33a}
        \mprob{\sum_{k=1}^t(\mC k -\mnC k) > N^{1/2-u_1}} \leq N^{-v_1}.
    \end{align}
\end{prop}
We start with an outline of the proof Proposition \ref{coupling of NR:B}.

The difference between $\mC k$ and $\mnC k$ stems from individuals, whose mark has appeared
previously and which is consequently discarded together with all its descendants. Call such an individual a \emph{duplicate}.
To show the claim \eqref{mainprop33a} we need some more details
of the thinning procedure of the \mRNprocess\ $\{\mZ t, \mJ t\}_{t\geq 0}$.
We relabel the marks of the \mRNprocess\
 $$
 \mM[0] 1,\mM [1] 1 ,\ldots, \mM [1]{\mZ 1},\mM [2]{1},\ldots,
 $$
given in Section \ref{subsec:static}, as
\begin{align} \label{trnaslate}
 M_0, M_1,M_2 ,\ldots, M_{\mZ 1},M_{\mZ 1 + 1},\ldots\
 \end{align}
By definition, $M_0$ is a random variable uniformly chosen from $\{1,2,\ldots,N\}$ and the marks $M_v$, for $v> 0$,
are independent copies of the random variable $M$ given by \eqref{dist mark}.
Introduce the random variable $Y_v$,
such that $Y_v=1$ if $M_v$ is a duplicate and $Y_v=0$ otherwise, so that
\begin{align} \label{def Y}
  Y_v = \mindic{\cup_{w=0}^{v-1} \{  M_v = M_w\}},
  \end{align}
and denote by $s(v)$, $v\geq 0$, the generation of individual $v$ in the \mRNprocess. Let $\text{Dup}_t$ be the set of all the duplicates in the first $t$ generations of the \mRNprocess, so that
\begin{equation}
\label{def_dup}
\text{Dup}_t =\bra{v\geq 1: Y_v =1 \text{ and } s(v) \leq t}, \,\, t \geq 0.
\end{equation}

Consider the subtree with root the duplicate $d$,
$d\in\text{Dup}_t$. The \mRNprocess\ is a marked BP, therefore the
subtree with root $d$ is also a marked BP and we denote this
subtree by $\{Z^{\smallsup{d}}_t, \mJd d t\}_{t\geq 0}$.
\ch {Observe that the mark of a duplicate, which we denote by the random variable $M_d|Y_d=1,$
depends on the marks seen previously; the distribution of the mark of a duplicate is biased, because of the duplication.
As a consequence also the offspring distribution of the first generation  of $\{Z^{\smallsup{d}}_t\}_{t\geq 0}$
differs  and we will denote this distribution by
    \begin{align} \label{def:off_dup}
        \hat f_n^{\smallsup{d}}
                =\mprob{\mPoisson{\mli {M_d}}=n|Y_d=1}, \,\, n\geq 0,
    \end{align}
  where $\lambda(a)=\lambda_a$ for any $a$.
In second and further generations the marks are independent copies
of the mark $M$ with distribution given by \eqref{dist mark} and
the offspring are independent with common distribution $\{\mg n\}_{n\geq 0}$, defined in \eqref{bepaling g}.
Therefore, denoting by
$$
    \mJd dk=(\mMd dk1,\mMd dk2,\ldots,\mMd dk{\mZd dk}), \quad k\geq 0,
$$
the vector of marks of the individuals in generation $k$ of the subtree with root $d$, we have for $k=0$, that $\mMd d01$ has distribution
$M_d|Y_d=1,$ whereas for $k\ge 1$, all the marks $\mMd dkj$ are independent with common distribution given by the right side of \eqref{dist mark}.}

By construction,
    \begin{align}\label{eq:verschilZ}
    \sum_{k=0}^t(\mZ k -\mnZ k)
            &\leq \sum_{d\in\text{Dup}_t}\sum_{k=0}^{t-s(d)}Z^{\smallsup{d}}_k
            = |\text{Dup}_t| + \sum_{d\in\text{Dup}_t}\sum_{k=1}^{t-s(d)}Z^{\smallsup{d}}_k,&
\intertext{and, similarly, } \label{eq:verschilC}
    \sum_{k=0}^t(\mC k -\mnC k)&
    \leq \sum_{d\in\text{Dup}_t}
            \lambda(\mMd d0 1)
                + \sum_{d\in\text{Dup}_t}\sum_{k=1}^{t-s(d)}\sum_{m=1}^{\mZd dk} \mli{\mMd dk m}.&
\end{align}

Suppose \ch{that} we can find an event $\cal D$, with $\mprob{{\cal D}^c}=\mO{ N^{-\eta/2}}$, such that the event $\cal D$
ensures that there are \ch {only a} few mismatches between the \mRNprocess\ and \mrRNprocess.
Then the Markov inequality, yields:
    \begin{align}\label{doel pre:eq}
    \mprob{\sum_{k=1}^t(\mC k -\mnC k) > N^{1/2-u_1}}
    \leq N^{-1/2+u_1}\mexpec{\sum_{k=1}^t(\mC k -\mnC k)\mindic{\cal D}} + \mO{ N^{-\eta/2}}.
    \end{align}
According to (\ref{doel pre:eq}) it is hence sufficient to show that
\begin{align} \label{eq:goal:dups}
    \mexpec{\sum_{k=1}^t(\mC k -\mnC k)\mindic{\cal D}}
            =\mO{N^{1/2-u_1-v_1}}.
\end{align}

The capacity of each duplicate is bounded from above  by $N^\gamma$, as a consequence of (\mCo3).
Furthermore, as the marks in the second and further
generations of a duplicate are independent copies of the random variable $M$,
and $\mexpec{\ml M}=\nu_\sN$, the expected value of \eqref{eq:goal:dups} becomes
\begin{align} \label{doel:eq}
    \mexpec{\sum_{k=1}^t(\mC k -\mnC k)\mindic{\cal D}}
        &\leq N^{\gamma}\mexpec{|\text{Dup}_t|\mindic{\cal D}}
+ \nu_\sN \mexpec{\sum_{d\in\text{Dup}_t}\sum_{k=1}^{t-s(d)}Z^{\smallsup{d}}_k\mindic{\cal D}}.&
\end{align}
Using auxiliary lemmas we will show that
\begin{align}\label{doel:sub:eq}
&\mexpec{|\text{Dup}_t|\mindic{\cal D}}
        = \mO{N^{1/2-\rho-\gamma}}&
&\text{and}&
&\mexpec{\sum_{d\in\text{Dup}_t}\sum_{k=1}^{t-s(d)}Z^{\smallsup{d}}_k\mindic{\cal D}}
        =\mO{N^{1/2-\rho}},&
\end{align}
for some $\rho>0$. Combining \eqref{doel:eq} and \eqref{doel:sub:eq} yields \eqref{eq:goal:dups}.
We end the outline with a list of all statements that we intend to prove and which
together imply the statement in Proposition \ref{coupling of NR:B}. The list consists of the following steps:
    \begin{enumerate}
        \item Define $\cal D$ and show that $\mprob{{\cal D}^c}\le N^{-\eta/2}$,
        \item Show the two statements in (\ref{doel:sub:eq}).
    \end{enumerate}

To control the size of the set
$\text{Dup}_t$ we will use  two lemmas.
\begin{lemma} \label{LA.2.4}
    For $\eta,\delta\in(-1/2,1/2)$ and all $t \leq (1/2+\eta)\log_\nu N$,
    $$
        \mprob{\sum_{k=0}^t \mZ k > N^{1/2+\delta}}=
        \mO{(\log_\nu N)N^{-(\delta-\eta)}}.
        $$
Consequently,
if $\{\SPT_t\}_{t\geq 0}$ are the reachable sets of an uniformly chosen node $A$ in the graph $G_\sN$, then
$$
    \mprob{|\SPT_t|> N^{1/2+\delta}}
    =
    \mO{(\log_\nu N)N^{-(\delta-\eta)}} .
    $$
\end{lemma}

\noindent{\bf Proof.}
From  (\mCo1), and for $k \leq c \log_\nu N$,
\eqn{ \label{cons:nu:afs}
\nu_\sN^k =\mO{N^c},
}
Taking $c=\tfrac12+\eta,$ and using \eqref{reqA} and \eqref{cons:nu:afs}, we find for $t \leq (1/2+\eta) \log_\nu N$,
\begin{align}
    \mprob{\sum_{k=0}^t\mZ k > N^{1/2+\delta}}
&\leq  N^{-1/2-\delta}\sum_{k=0}^t\mexpec{\mZ k }
        =N^{-1/2-\delta}\sum_{k=0}^t\mu_\sN \nu_\sN^k
       \nonumber \\&=N^{-1/2-\delta}\sum_{k=0}^t \mO{N^{1/2+\eta}}
        =\mO{(\log_\nu N) N^{\eta-\delta}}.
\label{verg2}
\end{align}
Using the coupling \eqref{kop nr ne},
\begin{align*}
        \mprob{|\SPT_t|> N^{1/2+\delta}}
            &=\mprob{\sum_{k=0}^t \mnZ k > N^{1/2+\delta}}
            \leq \mprob{\sum_{k=0}^t \mZ k > N^{1/2+\delta}}
            =\mO{(\log_\nu N)N^{\eta-\delta}}.
\end{align*}
\qed

\begin{lemma} \label{NO DUP} For each $\delta \in (-1/2,1/2)$, $u \geq 0$ and $v \leq N^{1/2+\delta}$,
\begin{align} \label{dup0}
        \mprob{{{\sum_{w=1}^v Y_w \geq N^{u}}}}
                = \mO{N^{-(u - 2\delta) }}.
\end{align}
\end{lemma}

\noindent{\bf Proof.}
From the definition of $Y_v$, see \eqref{def Y}, and Boole's inequality,
\begin{align*}
    \mprob{Y_v = 1|M_v=m}
\nonumber&
            = \mprob{\cup_{w=0}^{v-1} \bra{M_w=m}}
\leq \sum_{w=0}^{v-1} \mprob{M_w=m}
            = \frac{1}{N} + (v-1)\frac{\ml m}{l_\sN}.
\end{align*}
Therefore, using $v\leq N^{1/2+\delta}$, we have that
\begin{align*}
\mprob{Y_v = 1}&
        =\frac{1}{N}\sum_{m=1}^N  \frac{\ml m}{l_\sN} + (v-1) \sum_{m=1}^N   \frac{\ml m^2}{l_\sN^2}
        =\frac{1}{N}+ \frac{v-1}{N}\frac{\mu_\sN}{\nu_\sN}
        \leq \frac{1}{N}+ \frac{\mu_\sN}{\nu_\sN} N^{\delta-1/2}
        = \mO{ N^{\delta-1/2}}.
\end{align*}
From the Markov inequality, \ch{and since $v\leq N^{1/2+\delta}$, we find}
\begin{align} \nonumber
        \mprob{\sum_{w=1}^v Y_w \geq N^{u}}
        &
        \leq N^{-u} \sum_{w=1}^v \mprob{Y_w=1}
        = \mO{v  N^{\delta-1/2-u}}=\mO{N^{2\delta-u}}. \ch{\qquad \mbox{\qed}}
\end{align}

\bigskip

\noindent{\bf Proof of Proposition \ref{coupling of NR:B}.}
Define the event ${\cal D}={\cal D}_1 \cap {\cal D}_2 \cap {\cal D}_3$, where
\begin{align}
\label{def-D}
    {\cal D}_1&=\bra{  \text{Dup}_{\bfloor{(1/2-2\eta)\log_\nu N}} =    \emptyset},&
    \\
    {\cal D}_2&=\bra{\brc||{\text{Dup}_{\bfloor{(1/2+\eta)\log_\nu N}}}     \leq N^{5\eta}},&
    \\
    {\cal D}_3&=\bigcap_{d\in\text{Dup}_{\bfloor{(1/2+\eta)\log_\nu N}}}\bra{ Z^\smallsup{d}_1 < N^{\theta}},
\end{align}
and where $\theta=6\eta+\gamma$.

It remains to prove that $\mprob{{\cal D}^c}=\mO{ N^{-\eta/2}}$ (Step 1), with $\cal D$ as defined
above, and the two statements in (\ref{doel:sub:eq}) (Step 2).

We start with the first statement in (\ref{doel:sub:eq}). We have,
\begin{eqnarray}
    \label{firstA312}
    \mexpec{|\text{Dup}_{t}|\mindic{\cal D}}
        \leq    \mexpec{|\text{Dup}_{\bfloor{(1/2+\eta)\log_\nu N}}|\mindic{\cal D}}
        \leq N^{5\eta}\leq N^{1/2-\rho-\gamma},
\end{eqnarray}
by choosing $\eta$ and $\rho$ small and  the fact that $\gamma\in(1/(\tau-1), 1/2)$, c.f. \eqref{P:assumP}.

The second statement in (\ref{doel:sub:eq}) proceeds in the following way.
On the event ${\cal D}_1 \cap {\cal D}_2$ duplicates do not appear
in the first $\lfloor (1/2-2\eta)\log_\nu N\rfloor$ generations, implying
that on ${\cal D}_1 \cap {\cal D}_2$,
\begin{align}\label{maxdepdup}
\max_{d\in\text{Dup}_{t}} (t-s(d))
       \leq t - \lfloor (1/2-2\eta)\log_\nu N\rfloor
                \leq 3\eta \log_\nu N,
\end{align}
since $t\leq (1/2+\eta)\log_\nu N$.
Thus,
\begin{equation}
\label{secondA312}
\mexpec{\sum_{d\in\text{Dup}_{t}}\sum_{k=1}^{{t}-s(d)}Z^{\smallsup{d}}_k\mindic{\cal D}}
\leq
\mexpec{\sum_{d\in\text{Dup}_{t}}\sum_{k=1}^{\bfloor{3\eta \log_\nu N}}Z^{\smallsup{d}}_k\mindic{\cal D}}.
\end{equation}
The right side is the expected size of the progeny of the duplicates and their offspring.
The expectation of the total number of children of all the duplicates on ${\cal D}_3$ is bounded from above by the
expectation of the total number of duplicates times $N^\theta$.
Since on ${\cal D}_2$ the total \ch{expected } number of duplicates is bounded from above by $N^{5\eta}$,
the total number of children of all the duplicates is bounded from above by $N^{5\eta+\theta}$ on ${\cal D}_2 \cap {\cal D}_3$.
Since, the mark of a child of a duplicate is by definition an independent copy of the random variable $M$,
the offspring distribution of each child of a duplicate is an independent copy of $\{\hat Z_l\}_{l\geq 0}$,
where $\{\hat Z_l\}_{l\geq 0}$ is a BP with $\hat Z_0=1$ and where each individual has offspring distribution $\{\mg n\}_{n\geq 0}$.
Therefore, the right side of \eqref{secondA312} can be bounded by
$$
\mexpec{\sum_{d\in\text{Dup}_{t}}\sum_{k=1}^{{t}-s(d)}Z^{\smallsup{d}}_k\mindic{\cal D}}
        \leq
N^{5\eta+\theta}+N^{5\eta+\theta}\mexpec{\sum_{k=1}^{\bfloor{3\eta \log_\nu N}}\hat Z_k}.
$$
Using \eqref{cons:nu:afs},
$$
\mexpec{\sum_{k=1}^{\bfloor{3\eta (\log_\nu N)}}\hat Z_k}
        =\sum_{k=1}^{\bfloor{3\eta (\log_\nu N)}}\nu_\sN^k
        =\sum_{k=1}^{\bfloor{3\eta (\log_\nu N)}} \mO{N^{3\eta}}
        =\mO{(\log_\nu N )N^{3\eta}}=\mO{N^{4\eta}}.
$$
Thus, we can bound \eqref{secondA312} by
$$
\mexpec{\sum_{d\in\text{Dup}_{t}}\sum_{k=1}^{{t}-s(d)}Z^{\smallsup{d}}_k\mindic{\cal D}}
\leq
N^{5\eta+\theta}+N^{5\eta+\theta}\mO{N^{4\eta}}=\mO{N^{\theta+9\eta}}=\mO{N^{15\eta+\gamma}},
$$
which is sufficient for the second statement of (\ref{doel:sub:eq}), since $\gamma < 1/2$ and we can pick $\eta$
arbitrarily small. This completes the proof of Step 2.

We continue with the proof of Step 1. For this, we bound
    \eq
    \label{splitprob}
    \prob({\cal D}^c)\leq \prob({\cal D}_1^c)+\prob({\cal D}_2^c)+\prob({\cal D}_3^c \cap {\cal D}_2).
    \en
We now bound each of these terms one by one. For $\prob({\cal D}_1^c)$, we use Lemma
\ref{LA.2.4} and Lemma \ref{NO DUP} to obtain
\begin{align}
\mprob{{\cal D}^c_1}&
        \leq\mprob{{\cal D}^c_1 \cap \bra{\sum_{k=0}^{\bfloor{(1/2-2\eta)\log_\nu N}} \mZ k < N^{1/2-\eta}}}
        +\mO{(\log_\nu N)N^{-\eta}}\nonumber
\\&
        \leq \mprob{\sum_{w=1}^{\bfloor{N^{1/2-\eta}}}Y_w > 0}+\mO{N^{-\eta/2}}
= \mO{N^{-\eta/2}},
\label{boundD1}
\end{align}
and, similarly,
\begin{align}
    \mprob{{\cal D}^c_2}&
            \leq \mprob{{\cal D}^c_2 \cap       \bra{\sum_{k=0}^{\bfloor{(1/2+\eta)\log_\nu N}} \mZ k < N^{1/2+2\eta}}}
            +\mO{(\log_\nu N)N^{-\eta}}\nonumber
\\&
        \leq \mprob{\sum_{w=1}^{\bfloor{N^{1/2+2\eta}}}Y_w > N^{5\eta}}+\mO{N^{-\eta/2}}
        = \mO{N^{-\eta/2}}.
\label{boundD2}
\end{align}

Using the Markov inequality and Boole's inequality
we have that
\begin{align} \label{adsfa}
\mprob{{\cal D}^c_3 \cap {\cal D}_2} \leq
N^{-\theta}\mexpec{\sum_{d\in\text{Dup}_{k^*}}Z^\smallsup{d}_1 \mindic{{\cal D}_2}},
\end{align}
where, for convenience, we set $k^*=\bfloor{(1/2+\eta)\log_\nu N}$.
If we condition on the sequence $\{M_v\}_{v\geq0}$ then $Z^\smallsup{d}_1 $, for $d\in \text{Dup}_{k^*}$,
is an independent Poisson random variable with mean $\mli{M_d}$.
Thus, conditioned on the sequence $\{M_v\}_{v\geq0}$, the sum $\sum_{d\in\text{Dup}_{k^*}}Z^\smallsup{d}_1$
is distributed as a Poisson random variable with mean $\sum_{d\in\text{Dup}_{k^*}}M_d$.
In turn this sum can be stochastically bounded from above by a Poisson random variable with mean $|\text{Dup}_{k^*}|N^{\gamma}$,
because $\mli{M_d} \leq \ml \sN ^\smallsup{\sN} \leq N^{\gamma}$ for all $d\in\text{Dup}_{k^*}$ by (\mCo3).
Since ${\cal D}_2$ and $\text{Dup}_{k^*}$ are measurable with respect to $\{M_v\}_{v\geq0}$, the above yields
\begin{align}
\mexpec{\sum_{d\in\text{Dup}_{k^*}}Z^\smallsup{d}_1 \mindic{{\cal D}_2}}
&\nonumber
        = \mexpec{\sum_{d\in\text{Dup}_{k^*}}\mexpec{\brc.|{\mPoisson{Z^\smallsup{d}_1}}\{M_v\}_{v\geq0} }\mindic{{\cal D}_2}}
    \\\nonumber&= \mexpec{\mexpec{\brc.|{\mPoisson{\sum_{d\in\text{Dup}_{k^*}}Z^\smallsup{d}_1}}\{M_v\}_{v\geq0} }\mindic{{\cal D}_2}}\\
        &\nonumber\leq \mexpec{\mexpec{\brc.|{\mPoisson{|\text{Dup}_{k^*}| N^\gamma}}\{M_v\}_{v\geq0} }\mindic{{\cal D}_2}}
        \\&\leq \mexpec{\mexpec{\brc.|{\mPoisson{ N^{5\eta+\gamma}}}\{M_v\}_{v\geq0} }}
        =\mexpec{\mPoisson{ N^{5\eta+\gamma}}}
            =N^{5\eta+\gamma}.
\end{align}
Thus, we bound \eqref{adsfa} by
\begin{align} \label{boundD3}
    \mprob{{\cal D}^c_3 \cap {\cal D}_2}
            \leq N^{-\theta}\mexpec{\sum_{d\in\text{Dup}_{k^*}}Z^\smallsup{d}_1 \mindic{{\cal D}_2}}
            \leq  N^{-\theta+5\eta+\gamma}=N^{-\eta}.
\end{align}
Hence, by substituting (\ref{boundD1}), (\ref{boundD2}) and (\ref{boundD3})  in (\ref{splitprob}) we obtain,
$$
\mprob{{\cal D}^c}=\mO{N^{-\eta/2}},
$$
which proves Step 1 and therefore the proposition.
\qed

\subsection{Proof of Proposition \ref{coupling of sums}} \label{proposition:coupling of sums:bewijs}
In this section we prove Proposition \ref{coupling of sums}, which we restate here for convenience as
Proposition \ref{B:coupling of sums}.

\begin{prop} \label{B:coupling of sums}
There exist constants $\alpha_3,\beta_3,\eta>0$ such that for all $t \le (1+2\eta)\log_{\nu} N$,
as $N\to\infty$,
\begin{eqnarray} \label{coupling of sums1:eq:B}
\prob\Big(\frac{1}{N}\Big|
\sum_{k=2}^{t+1} \msZi 1 {\lceil k/2 \rceil} \msZi 2 {\lfloor k/2 \rfloor}
-
\sum_{k=2}^{t+1}\mnCi 1 {\lceil k/2 \rceil} \mnCi 2 {\lfloor k/2\rfloor}\Big|>N^{-\alpha_3}\Big)=\mO{N^{-\beta_3}}.
\end{eqnarray}
\end{prop}
We start with an outline of the proof.
Define
\begin{align}
\label{a44}
T_1&=\frac{1}{N} \sum_{k=2}^{t+1} {\msZi 2 {\bfloor{k/2}} \brc||{\msZi 1 {\bceil{k/2}}-\mnCi 1 {\bceil{k/2}}}}&
&\text{and}&
T_2&=\frac{1}{N} \sum_{k=2}^{t+1} {\mnCi 1 {\bceil{k/2}} \brc||{\msZi 2 {\bfloor{k/2}}-\mnCi 2 {\bfloor{k/2}}}}.&
\end{align}
To show \eqref{coupling of sums1:eq:B}, it suffices to prove that
for an appropriate chosen event $\cal H$,
\begin{align} \label{1:to proof}
        &\mprob{T_i\mindic{\cal H} >\frac{1}{2}N^{-\alpha_3}}=\mO{N^{-\beta_3}},\qquad
i=1,2,
\end{align}
and that $\mprob{{\cal H}^c}= \mO{N^{-\xi}}$ for some $\xi>0$.
We choose the event $\cal H$ such that on this event, for each $k \leq t \leq  (1+2\eta)\log_\nu N$,
\begin{align} \label{as1}
    \mexpec{\msZi 2 {\bfloor{k/2}}\mindic{\cal H}}&=\mO{N^{1/2+2\eta}},&
&&
    \mexpec{\brc||{\msZi 1 {\bceil{k/2}} - \mnCi 1 {\bceil{k/2}}}\mindic{\cal H}}&=\mO{N^{1/2-u}},&
\\
\label{as2}
\mexpec{\mnCi 1 {\bceil{k/2}}\mindic{\cal H}}&=\mO{N^{1/2+2\eta}}&
&\text{and}&
\mexpec{\brc||{\msZi 2 {\bfloor{k/2}} - \mnCi 2 {\bfloor{k/2}}}\mindic{\cal H}}&=\mO{N^{1/2-u}},&
\end{align}
for some $u,\alpha_3,\eta>0$ such that
\begin{align} \label{cond finale}
\alpha_3 < u -2\eta.
\end{align}
The claims \eqref{as1}, \eqref{as2} and \eqref{cond finale} imply \eqref{1:to proof}.
This is straightforward from the Markov inequality, the independence of $\msZi 1 \nill$ between $\msZi 2\nill$ and
the associated capacities and the inequality $t \leq  (1+2\eta)\log_\nu N$. We leave it to the reader to verify this.

Instead of showing \eqref{as1} and \eqref{as2}, it is sufficient to show that there exist constants
$\alpha_3,\eta,u>0$ satisfying \eqref{cond finale} and such that for each $k \leq t\leq (1/2+\eta)\log_\nu N$,
\begin{align} \label{as}
        \mexpec{\brc||{\msZ k - \mnC k}\mindic{\cal H}}&=\mO{N^{1/2-u}},&
        \mexpec{\msZ k\mindic{\cal H}}&=\mO{N^{1/2+2\eta}}&
&\text{and}&
        \mexpec{\mnC k\mindic{\cal H}}&=\mO{N^{1/2+2\eta}},&
\end{align}
which simplifies the notation considerably.

We choose  ${\cal H}={\cal H}_1 \cap {\cal H}_2 \cap {\cal H}_3$, where
\begin{align} \label{def F}
        {\cal H}_1 &=  \bra{\sum_{k=1}^t(\mC k -\mnC k) \leq N^{1/2-u_1}},&\\
        {\cal H}_{2,k} &= \bra{\brc||{\mC k - \nu_\sN\mZ{k-1}}< N^{1/2-\delta}},& {\cal H}_2 &= \cap_{k=1}^t{\cal H}_{2,k},\\
        {\cal H}_{3,k} &= \bra{\brc||{\msZ k-\nu \msZ{k-1}}< N^{1/2-\delta}},& {\cal H}_3       &= \cap_{k=1}^t{\cal H}_{3,k}.
\end{align}
Hence, for the proof of Proposition \ref{B:coupling of sums} it suffices to show the three claims in \eqref{as} together with $u-2\eta>0$, and
that for some $\xi>0$ we have that $\mprob{{\cal H}^c}=\mO{N^{-\xi}}$.
For the latter statement, we use the following lemmas.
In the \ch{first} lemma we denote by $d_\sN$ a random variable
with distribution $\{\mg n\}_{n\geq0}$ given by \eqref{bepaling g}.


\begin{lemma}[Bound the variance of the marked capacity]
\label{LA.2.5}
If the conditions (\mCo1)--(\mCo3) do hold then
\begin{align} \label{LA.2.5e}
    \mvar{\ml{M}}& =\mO{N^{\gamma}}&
    &\text{and}&
    \mvar{d_\sN} &=\mO{N^{\gamma}},&
\end{align}
where $\gamma$ was introduced in \eqref{P:assumP}.
\end{lemma}

\begin{lemma}[Strong convergence of $\{\mg n\}_{n\geq0}$ to $\{\mgl n\}_{n\geq0}$]
\label{lem-ngnNgn}
Assume that conditions given by (\mCo1)--(\mCo3) do hold. Then,
there exists an $\alpha>0$ such that
    \eqn{
    \sum_{n=0}^{\infty} n |\mg n-\mgl n|
    =O(N^{-\alpha}).
    }
\end{lemma}
The proofs of the above lemmas are deferred to the end of this section. We start with the proof
of Proposition \ref{B:coupling of sums}, given the Lemmas \ch{\ref{LA.2.5}-\ref{lem-ngnNgn}.}

\bigskip

\noindent{\bf Proof of Proposition \ref{B:coupling of sums}.}
We have to prove that the three claims of \eqref{as} hold, with $u -2\eta>0$.
Consider the first claim given by \eqref{as}. Using the triangle inequality we arrive at
\begin{multline}
\label{grote afschatter}
\mexpec{\brc||{\mnC k-\msZ k}\mindic{\cal H}}
    \leq
    \mexpec{\brc||{\mnC k-\mC k}\mindic{{\cal H}}}
    +
    \mexpec{\brc||{\mC k -\nu_\sN \mZ{k-1}} \mindic{{\cal H}}}
    \\+
    \mexpec{\brc||{\mZ{k-1}(\nu_\sN-\nu)} \mindic{{\cal H}}}
    +
    \nu\mexpec{\brc||{\mZ{k-1} - \msZ{k-1}}\mindic{{\cal H}}}
    +
    \mexpec{\brc||{\nu\msZ{k-1} - \msZ{k}}\mindic{{\cal H}}}.
\end{multline}
The first, second and the last term on the right side of \eqref{grote afschatter}, we bound by $N^{1/2-\min{\{\delta,u_1\}}}$,
using the events ${\cal H}_1$, ${\cal H}_2$  and ${\cal H}_3$, respectively.
We bound the third term of \eqref{grote afschatter} using \eqref{cons:nu:afs}, which gives
\begin{align}\label{afs zoveel}
\mexpec{\brc||{\mZ{k-1}(\nu_\sN-\nu)} }
\leq N^{-\alpha_1}\mexpec{\mZ{k-1}}
=N^{-\alpha_1}\nu_\sN^{k-1}
=\mO{N^{1/2+\eta-\delta}}.
\end{align}
Finally, we need to bound the fourth term of \eqref{grote afschatter}.
The following inequality holds for each $k \leq t$ (compare with \cite[(A.1.4) and (A.1.15)]{finstub}),
\begin{align}\label{msZ i - mZ i}
\mexpec{\brc||{\mZ k-\msZ k}}\leq \max\bra{\nu-\alpha_\sN,\nu_\sN-\alpha_\sN}
\sum_{l=1}^t\mexpec{\mZ l}\br{\max\bra{\nu,\nu_\sN}}^{t-l},
\end{align}
where
\begin{align*}
\alpha_\sN=\sum_{n=0}^\infty n \min\{\mgl n, \mg n\}.
\end{align*}
We bound the sum in \eqref{msZ i - mZ i}, using \eqref{reqA} and \eqref{cons:nu:afs}:
\begin{align}  \nonumber
\sum_{l=1}^t \mexpec{\mZ l}\br{\max\bra{\nu,\nu_\sN}}^{t-l}
&=\sum_{l=1}^t \mu_\sN\nu_\sN^{l-1}\br{\max\bra{\nu,\nu_\sN}}^{t-l}
\leq t {\max\bra{\mu,\mu_\sN}} \br{\max\bra{\nu,\nu_\sN}}^{t-1}
\\\label{Afs:1} &\leq t \br{1+\mO{N^{-\alpha_1}}}^{t} \mO{N^{1/2+\eta}} \leq  N^{1/2+2\eta},
\end{align}
for $N$ sufficiently large.
Furthermore from Lemma \ref{lem-ngnNgn},
\begin{align} \label{Afs:2}
\max\bra{\nu-\alpha_\sN,\nu_\sN-\alpha_\sN}&\leq \sum_{n=1}^\infty n \brc||{\mgl n-\mg n}
\leq N^{-\alpha_1},
\end{align}
for some $\alpha_1>0$.
Combining \eqref{msZ i - mZ i}, \eqref{Afs:1} and \eqref{Afs:2}, we obtain that,
\begin{align*}
\mexpec{\brc||{\mZ k-\msZ k}}  =\mO{N^{1/2+2\eta-\alpha_1}}.
\end{align*}
All together, the left side of \eqref{grote afschatter} satisfies
\begin{align}\label{sdoel1}
\mexpec{\brc||{\mnC k-\msZ k}\mindic{\cal H}}=\mO{N^{1/2-\min{\{\delta,u_1,\delta-\eta,\alpha_1-2\eta\}}}}.
\end{align}
This \ch{yields} the first claim in \eqref{as} with $$u=\min{\{\delta,u_1,\delta-\eta,\alpha_1-2\eta\}}.$$

The second claim of \eqref{as} is evident, because $k\leq t \leq (1/2+\eta)\log_\nu N$ and therefore
\begin{align}\label{controZ}
\mexpec{\msZ k\mindic{\cal H}}\leq\mexpec{\msZ k}=\mu\nu^{k-1}\leq \mu\nu^{t-1} = \mO{N^{1/2+\eta}}.
\end{align}
Finally, the third claim of \eqref{as} follows from
using that $\mnC k\leq \mC k$,
\begin{align*}
\mexpec{\mnC k\mindic{\cal H}}\leq \mexpec{\mnC k}\leq \mexpec{\mC k}&=
\mexpec{\sum_{v=1}^{\mZ {k-1}} \mli{\mM[k-1] v}}.
\end{align*}
Now by taking conditional expectation with respect to $\mZ {k-1}$
and the capacities,  we obtain
\begin{align*}
\mexpec{\brc.|{\sum_{v=1}^{\mZ {k-1}} \mli{\mM[k-1] v}}\mZ {k-1}}
&=
\sum_{v=1}^{\mZ {k-1}}\mexpec{\brc.|{\mli{\mM[k-1] v}}\mZ {k-1}}
=\mZ {k-1}\mexpec{\ml M}
=\nu_\sN \mZ {k-1},
\end{align*}
so that \eqref{cons:nu:afs} implies that
\begin{align*}
\mexpec{\mnC k\mindic{\cal H}}&
\leq
\mexpec{\mexpec{\brc.|{\sum_{v=1}^{\mZ {k-1}} \mli{\mM[k-1] v}}\mZ {k-1}}}
= \mexpec{ \nu_\sN\mZ {k-1}}
=\mO{N^{1/2+\eta}}.
\end{align*}
Thus, we have shown the claims given by \eqref{as} if we restrict $\eta$ to
$0 <\eta < \min\{\delta,\alpha_1/2\}$ due to \eqref{sdoel1}.
To satisfy condition \eqref{cond finale}, we restrict $\eta$ to
\begin{align} \label{res eta}
0<\eta < \min\{\alpha_1/4,u_1/2,\delta/2\},
\end{align}
and  pick
$\alpha_3=(u-2\eta)/2$. Then $\alpha_3>0$ and condition \eqref{cond finale} is satisfied,
because
$$
u -\eta < \alpha_3=(u-2\eta)/2=\min\{\alpha_1-4\eta,
u_1-2\eta,\delta-2\eta\}/2,
$$
for each $0 < \delta < 1/2$, and for each $\eta$ given by \eqref{res eta}.

We finish the proof by showing that $\mprob{{\cal H}^c} =\mO{N^{-\xi}}$ for some $\xi>0$.
A simple bound is
\begin{align}
\label{Hcompl}
\mprob{{\cal H}^c}
\leq
\sum_{i=1}^3 \mprob{{\cal H}_i^c}.
\end{align}

We bound $\mprob{{\cal H}^c_1}$ using Proposition \ref{coupling of NR:B}, which states that there exist constants $u_1,v_1$ such that
\begin{align}\label{H1af}
\mprob{{\cal H}^c_1} =\mO{N^{-v_1}}.
\end{align}
Next, we will bound $\mprob{{\cal H}_2^c}$ and $\mprob{{\cal H}_3^c}$ from above.
\ch {Using Chebychev's inequality, we will show that}
\begin{align}
\mprob{{\cal H}_{2,k}^c } = \mO{N^{-1/2+2\delta+\gamma+2\eta}}\text{,\, for $k\leq t \leq (1/2+\eta)\log_\nu N$.}
\end{align}
Then by Boole's inequality we have,
\begin{align}
\mprob{{\cal H}_{2}^c} &\leq
\sum_{k=0}^{t}\mprob{{\cal H}_{2,k}^c}
=  \mO{t N^{-1/2+2\delta+\gamma+2\eta}}
=\mO{(\log_\nu N)N^{-(1/2-2\delta-\gamma-2\eta)}}.
\label{H2af}
 \end{align}

Fix $k$, $1\leq k \leq t\leq (1/2+\eta)\log_\nu N$, and denote by $V_v=\mli {\mM[k-1] v}-\nu_\sN$ for $v=1,2,\ldots,\mZ {k-1}$.
We have that $\mexpec{V_v}=0$, and that, conditionally on $\mZ {k-1}$,
the sequence $\{V_v\}_{v=1}^{\mZ {k-1}}$ is an i.i.d. sequence.
Hence, using Lemma \ref{LA.2.5},
\begin{align}
\mexpec{\brc.|{\br{\sum_{v=1}^{\mZ {k-1}} V_v}^{2}}\mZ {k-1}}\nonumber
&
=\sum_{v=1}^{\mZ {k-1}}\sum_{w=1}^{\mZ {k-1}}\mexpec{{ V_v V_w}}
 \nonumber
= \sum_{v=1}^{\mZ {k-1}}\sum_{w=1,w\not=v}^{\mZ {k-1}}\mexpec{ V_v}\mexpec{V_w}
+\sum_{v=1}^{\mZ {k-1}}\mexpec{{ V_v^2}}
\\&=\sum_{v=1}^{\mZ {k-1}}\mexpec{{ V_v^2}}
=\mZ {k-1}\mvar{\mli M} \label{A.4.23a}
=\mO{N^\gamma \mZ {k-1} }.
\end{align}
Therefore, and since by \eqref{cons:nu:afs}, $\mexpec{\mZ {k-1}}=\nu_\sN^{k-1}
=\mO{N^{1/2+\eta}}$,
\begin{align*}
\mvar{\sum_{v=1}^{\mZ {k-1}} V_v}&=\mexpec{\mexpec{\brc.|{\br{\sum_{v=1}^{\mZ {k-1}} V_v}^{2}}\mZ {k-1}}}
= \mO{N^\gamma\mexpec{\mZ {k-1}}}=
\mO{N^{1/2+\gamma+\eta}}.
\end{align*}
Thus, by the Chebyshev inequality,
\begin{align*}
\mprob{{\cal H}_{2,k}^c} &
=\mprob{\brc||{\sum_{v=1}^{\mZ {k-1}} V_v} \geq N^{1/2-\delta}}
\leq N^{2\delta-1}\mvar{\sum_{v=1}^{\mZ {k-1}} V_v }
=\mO{N^{-1/2+2\delta+\gamma+\eta}}.
\end{align*}

Similarly, we can show that
$$
\mprob{{\cal H}_{3,k}^c } =\mO{(\log_\nu N)N^{-(1/2-2\delta-\gamma-\eta)}},
$$
when we replace $\mZ {k-1}$ by $\msZ {k-1}$, set $V_n=X_{k-1,n}-\mexpec{d_\sN}$, where $X_{k-1,n}$ is an independent copy of $d_\sN$. Using Lemma \ref{LA.2.4} and \eqref{controZ}, this yields
\begin{align*}
\mprob{{\cal H}_{3,k}^c }
&=
\mprob{\brc||{\sum_{v=1}^{\msZ {k-1}} V_v} \geq N^{1/2-\delta}}
\leq
N^{2\delta-1}\mexpec{\mexpec{\brc.|{\br{\sum_{v=1}^{\msZ {k-1}} V_v}^{2}}\msZ {k-1}}}
\\&=N^{2\delta-1}\mexpec{\msZ {k-1}\mvar {d_\sN}}=N^{2\delta-1}\mvar {d_\sN}\mexpec{\msZ {k-1}}=\mO{N^{2\delta-1/2+\eta+\gamma}}.
\end{align*}
The above yields, compare with \eqref{H2af},
\begin{align}
\mprob{{\cal H}_{3}^c} &
=\mO{N^{-(1/2-2\delta-\gamma-2\eta)}}.
\label{H3af}
 \end{align}

Combining \eqref{H1af}, \eqref{H2af}, \eqref{H3af} with \eqref{Hcompl} gives
$$
\mprob{{\cal H}^c}=\mO{N^{-\xi}},
$$
where $\xi=\min\{(1/2-2\delta-\gamma-2\eta), v_1\}$.
Remember that $1/2 - \gamma>0$, that
$\delta$ is restricted to $0<\delta<1/2$ and $\eta$ is restricted by \eqref{res eta}.
So, pick $\eta,\delta>0$ such that  $2\delta+\eta < 1/2-\gamma$,
because we can pick $\delta>0$ and $\eta>0$ arbitrary small, then $\xi>0$.
\qed

\bigskip

\noindent{\bf Proof of Lemma \ref{LA.2.5}: }
Consider the first claim of \eqref{LA.2.5e}. The variance of a random variable is bounded from above by its second moment. Therefore, using \eqref{dist mark}, \eqref{reqA} and \eqref{reqC},
\begin{align*}
    {\mvar{\ml{M}}}&\leq {\mexpec{\ml{M}^2}}
    ={\frac{1}{l_\sN} \sum_{m=1}^N \ml{m}^3}
    \leq {\frac{\ml \sN ^\smallsup{\sN}}{l_\sN} \sum_{m=1}^N \ml{m}^2}
    \leq N^\gamma \nu_\sN
    = N^\gamma \br{\nu+\mO{N^{-\alpha_1}}}=\mO{N^\gamma}.
\end{align*}
We turn to the second claim of \eqref{LA.2.5e}. Using \eqref{helper:mom},
where we take $q=2$, and (\mCo3) we bound
\begin{align*}
\mvar{d_\sN}
        & \leq \sum_{n=1}^\infty n^2 \mg n
        =c+\frac{c}{\mu_{\sN}}\frac 1N\sum_{i=1}^N \lambda_i^{3}
        \leq c+\frac{c\ml n^{\smallsup{N}}}{\mu_{\sN}}\frac 1N\sum_{i=1}^N \lambda_i^{2}
        =c+c\nu_\sN \ml \sN^{\smallsup{N}}=\mO{N^{\gamma}}.
\end{align*}
\qed

\bigskip

{\noindent \bf Proof of Lemma \ref{lem-ngnNgn}:} \ch{By Fatou's Lemma,
as well as the convergence of
$\lim_{N\rightarrow \infty} \mg n=\mgl n$ which follows from (\mCo2), we have}
    \eqn{
    \sum_{n=0}^{\infty} n^q \mgl n
    =\sum_{n=0}^{\infty} n^q \liminf_{N\rightarrow \infty} \mg n
    \leq \liminf_{N\rightarrow \infty}\sum_{n=0}^{\infty} n^q \mg n<
    \infty,
    }
for $q<\tau-2$ by Lemma \ref{lem-gvslambda} combined with (\mCo3).

Then, we split, for $\zeta>0$,
    \eqn{
    \label{boundonsum}
    \sum_{n=0}^{\infty} n |\mg n-\mgl n|
    \leq \sum_{n=0}^{\lfloor N^{\zeta}\rfloor} n |\mg n-\mgl n|
    +\sum_{n=\lfloor N^{\zeta}\rfloor+1}^{\infty} n |\mg n-\mgl n|.
    }
Using (\mCo2), we bound
    \eqn{
    \sum_{n=0}^{\lfloor N^{\zeta}\rfloor} n |\mg n-\mgl n|
    \leq \lfloor N^{\zeta}\rfloor \sum_{n=0}^{\infty} |\mg n-\mgl n|
    \leq N^{\zeta-\beta}.
    }
Using (\mCo3) as well as \eqref{boundonsum}, we further obtain for $\vep>0$,
    \eqn{
    \sum_{n=\lfloor N^{\zeta}\rfloor}^{\infty} n |\mg n-\mgl n|
    \leq N^{-\zeta(\tau-3-\vep)}
    \sum_{n=\lfloor N^{\zeta}\rfloor}^{\infty} n^{\tau-2-\vep} (\mg n+\mgl n)
    =O(N^{-\zeta(\tau-3-\vep)}).
    }
Fix $\zeta<\beta$ and $\vep\in (0,\tau-3),$ which is possible since
$\tau>3$, then the claim follows with $\alpha=\min\{\beta-\zeta, \zeta(\tau-3-\vep)\}.$
\qed
}

\subsection*{Acknowledgements}
The work of HvdE and RvdH was supported in part by Netherlands Organisation for
Scientific Research (NWO). We thank Michel Dekking for carefully reading the manuscript and many suggestions, which resulted in a more accessible paper.
We thank Chris Klaassen for shedding light on equation \eqref{L1-stat}. Finally, we thank Olaf Wittich for discussions on Lemma \ref{lem-conv-int}.


\begin{thebibliography}{10}

\bibitem{asmussen}
S.~Asmussen.
\newblock Some martingale methods in the limit theory of supercritical
  branching processes.
\newblock In {\em Branching processes}, pages 1--26. Advances in Probability
  and related Topics, 1978.

\bibitem{BeC}
E.A. Bender and E.R. Canfield.
\newblock The asymptotic number of labelled graphs with given degree sequences.
\newblock {\em Journal of Combinatorial Theory}, A24(3):296--307, 1978.

\bibitem{bollobasboek}
B.~Bollob{\'a}s.
\newblock {\em Random Graphs}.
\newblock Cambridge University Press, 2001.

\bibitem{BBCR}
B.~Bollob{\'a}s, C.~Borgs, J.T. Chayes, and O.~Riordan.
\newblock Directed scale-free graphs.
\newblock In {\em Proceedings of the fourteenth annual ACM-SIAM symposium on
  Discrete algorithms}, pages 132--139. Symposium on Discrete Algorithms, 2003.

\bibitem{phasetrans}
B.~Bollob{\'a}s, S.~Janson, and O.~Riordan.
\newblock The phase transition in inhomogeneous random graphs.
\newblock {\em Random Structures and Algorithms}, 31(1):3--122, 2007.

\bibitem{Britton}
T.~Britton, M.~Deijfen, and A.~Martin-L{\"o}f.
\newblock Generating simple random graphs with prescribed degree distribution.
\newblock {\em Journal of Statistical Physics}, 124(6):1377--1397, 2006.

\bibitem{CLaverage}
F.~Chung and L.~Lu.
\newblock The average distances in random graphs with expected degrees.
\newblock {\em PNAS}, 99(25):15879--15882, 2002.

\bibitem{CLconnected}
F.~Chung and L.~Lu.
\newblock Connected components in random graphs with given expected degree
  sequences.
\newblock {\em Annals of Combinatorics}, 6(2):125--145, 2002.

\bibitem{EVshort}
H.~van~den Esker, R.~van~der Hofstad, and G.~Hooghiemstra.
\newblock Universality for the distance in finite variance random graphs.
\newblock {\em preprint}.

\bibitem{infmean}
H.~van~den Esker, R.~van~der Hofstad, G.~Hooghiemstra, and D.~Znamenski.
\newblock Distances in random graphs with infinite mean degrees.
\newblock {\em Extremes}, 8:111--140, 2006.

\bibitem{fellerb}
W.~Feller.
\newblock {\em An Introduction to Probability Theory and Its Applications
  Volume II, 2nd edition}.
\newblock John Wiley and Sons, New York, 1971.

\bibitem{Gut1}
A.~Gut.
\newblock {\em Probability: a graduate course}.
\newblock Springer, Berlin, 2005.

\bibitem{finstub}
R.~van~der Hofstad, G.~Hooghiemstra, and P.~Van Mieghem.
\newblock Distances in random graphs with finite variance degrees.
\newblock {\em Random Structures and Algorithms}, 27(2):76--123, 2005.

\bibitem{infvar}
R.~van~der Hofstad, G.~Hooghiemstra, and D.~Znamenski.
\newblock Distances in random graphs with finite mean and infinite variance
  degrees.
\newblock {\em Electronic Journal of Probability}, 12:703--766, 2007.

\bibitem{Jans08a}
S.~Janson.
\newblock Asymptotic equivalence and contiguity of some random graphs.
\newblock 2008.

\bibitem{JLR}
S.~Janson, T.~Luczak, and A.~Ruci{\'n}ski.
\newblock {\em Random Graphs}.
\newblock John Wiley \& Sons, New York, 2000.

\bibitem{MR95}
M.~Molloy and B.~Reed.
\newblock A critical point for random graphs with a given degree sequence.
\newblock In {\em Proceedings of the sixth international seminar on Random
  graphs and probabilistic methods in combinatorics and computer science},
  pages 161--179, 1995.

\bibitem{NSW00}
M.E.J. Newman, S.H. Strogatz, and D.J. Watts.
\newblock Random graphs with arbitrary degree distribution and their
  application.
\newblock {\em Physical Review E}, 64:026118, 2001.

\bibitem{norros1}
I.~Norros and H.~Reittu.
\newblock On the power-law random graph model of massive data networks.
\newblock {\em Performance Evaluation}, 55(1-2):3--23, 2004.

\bibitem{norros3}
I.~Norros and H.~Reittu.
\newblock On a conditionally $\text{P}$oissonian graph process.
\newblock {\em Advances in Applied Probability}, 38(1):59--75, 2006.

\bibitem{thorisson}
H.~Thorisson.
\newblock {\em Coupling, Stationarity and Regeneration}.
\newblock Springer, Berlin, 2000.

\end{thebibliography}
\end{document}